\documentclass{amsart}
\begin{document}
\title{Differential eigenforms}
\author{Alexandru Buium}
\def \bT{{\bf T}}
\def \cI{{\mathcal I}}
\def \cJ{{\mathcal J}}
\def \ZN{\bZ[1/N,\zeta_N]}
\def \tA{\tilde{A}}
\def \o{\omega}
\def \tB{\tilde{B}}
\def \tC{\tilde{C}}
\def \alph{A}
\def \bet{B}
\def \bsigma{\bar{\sigma}}
\def \y{^{\infty}}
\def \Ra{\Rightarrow}
\def \uBS{\overline{BS}}
\def \lBS{\underline{BS}}
\def \lB{\underline{B}}
\def \<{\langle}
\def \>{\rangle}
\def \hL{\hat{L}}
\def \cU{\mathcal U}
\def \cF{\mathcal F}
\def \S{\Sigma}
\def \st{\stackrel}
\def \sd{Spec_{\d}\ }
\def \pd{Proj_{\d}\ }
\def \s{\sigma_2}
\def \i{\sigma_1}
\def \bs{\bigskip}
\def \cD{\mathcal D}
\def \cC{\mathcal C}
\def \cT{\mathcal T}
\def \cK{\mathcal K}
\def \cX{\mathcal X}
\def \sX{X_{set}}
\def \cY{\mathcal Y}
\def \cS{X}
\def \cR{\mathcal R}
\def \cE{\mathcal E}
\def \tcE{\tilde{\mathcal E}}
\def \cP{\mathcal P}
\def \cA{\mathcal A}
\def \cV{\mathcal V}
\def \cM{\mathcal M}
\def \cN{\mathcal N}
\def \tcM{\tilde{\mathcal M}}
\def \caS{\mathcal S}
\def \cG{\mathcal G}
\def \cB{\mathcal B}
\def \tG{\tilde{G}}
\def \cF{\mathcal F}
\def \h{\hat{\ }}
\def \hp{\hat{\ }}
\def \tS{\tilde{S}}
\def \tP{\tilde{P}}
\def \tA{\tilde{A}}
\def \tX{\tilde{X}}
\def \tcS{\tilde{X}}
\def \tT{\tilde{T}}
\def \tE{\tilde{E}}
\def \tV{\tilde{V}}
\def \tC{\tilde{C}}
\def \tI{\tilde{I}}
\def \tU{\tilde{U}}
\def \tG{\tilde{G}}
\def \tu{\tilde{u}}
\def \chu{\check{u}}
\def \tx{\tilde{x}}
\def \tL{\tilde{L}}
\def \tY{\tilde{Y}}
\def \d{\delta}
\def \e{\chi}
\def \bZ{{\bf Z}}
\def \bV{{\bf V}}
\def \bF{{\bf F}}
\def \bE{{\bf E}}
\def \bC{{\bf C}}
\def \bO{{\bf O}}
\def \bR{{\bf R}}
\def \bA{{\bf A}}
\def \bB{{\bf B}}
\def \cO{\mathcal O}
\def \ra{\rightarrow}
\def \bx{{\bf x}}
\def \f{{\bf f}}
\def \bX{{\bf X}}
\def \bH{{\bf H}}
\def \bS{{\bf S}}
\def \bF{{\bf F}}
\def \bN{{\bf N}}
\def \bK{{\bf K}}
\def \bE{{\bf E}}
\def \bB{{\bf B}}
\def \bQ{{\bf Q}}
\def \bd{{\bf d}}
\def \bY{{\bf Y}}
\def \bU{{\bf U}}
\def \bL{{\bf L}}
\def \bQ{{\bf Q}}
\def \bP{{\bf P}}
\def \bR{{\bf R}}
\def \bC{{\bf C}}
\def \bD{{\bf D}}
\def \bM{{\bf M}}
\def \bP{{\bf P}}
\def \xtoleqr{x^{(\leq r)}}
\def \hU{\hat{U}}

\newtheorem{THM}{{\!}}[section]
\newtheorem{THMX}{{\!}}
\renewcommand{\theTHMX}{}
\newtheorem{theorem}{Theorem}[section]
\newtheorem{corollary}[theorem]{Corollary}
\newtheorem{lemma}[theorem]{Lemma}
\newtheorem{proposition}[theorem]{Proposition}
\theoremstyle{definition}
\newtheorem{definition}[theorem]{Definition}
\theoremstyle{remark}
\newtheorem{remark}[theorem]{Remark}
\newtheorem{example}[theorem]{\bf Example}
\numberwithin{equation}{section}
\address{University of New Mexico \\ Albuquerque, NM 87131}
\email{buium@math.unm.edu} \subjclass[2000]{11 F 32, 11 F 85}
\maketitle

\begin{abstract}
The aim of this paper is to show how  $\d-$characters of Abelian
varieties (in the sense of \cite{char}) can be used to construct
$\d-$modular forms of weight $0$ and order $2$ (in the sense of
\cite{difmod}) which are eigenvectors of Hecke operators. These
$\d-$modular forms have ``essentially the same" eigenvalues as
certain classical complex eigenforms of weight $2$.
\end{abstract}

\section{Introduction}
 The concept of $\d-${\it modular form} was
introduced in \cite{difmod}. Very roughly speaking a level one
$\d-$modular form of order $r$ is a ``homogeneous'' function of
plane  elliptic curves $y^2=x^3+ax+b$ (where $a,b \in
R:=\hat{\bZ}^{ur}_p$) that can be written as a $p-$adic restricted
power series in $a,b,\d a, \d b,...,\d^r a, \d^r b, \Delta^{-1}$,
where $\Delta:=4a^3+27b^2$ and  $\d^ia,\d^i b$ are the iterated
``Fermat quotients'' of $a,b$ with respect to $p$. We recall that
$\d x:=(\phi(x)-x^p)/p$, where $\phi:R \ra R$ is the lift of the
$p-$power Frobenius on $R/pR$. Morally one may view $\d$ as an
arithmetic analogue of a derivation (acting on ``numbers" rather
than ``functions") and one may view $\d-$modular forms as
``non-linear arithmetic differential operators of order $r$"
acting on pairs $(a,b)$. We shall review this concept presently,
from a slightly different (but equivalent) viewpoint.  There is a
level
 $N$ generalization of this.
Also there are Hecke operators $T(l)$ acting on $\d-$modular forms
(where $l$ are primes with $(l,Np)=1$) so one can talk about {\it
$\d-$ eigenforms} (for all these $T(l)$'s). Finally one can attach,
to $\d-$modular forms of order $r$,
$\d-${\it Fourier expansions} which are series
in the variables $q,q',...,q^{(r)}$.
 For applications of our theory we
refer to \cite{difmod}, \cite{book}.

There is an ``easy" way to construct $\d-$eigenforms  by
considering   $I-$linear combinations of ``$\phi-$powers'',
$f^{\phi^j}$, of  classical (complex) eigenforms $f$ where $I$ is
the ring generated by the isogeny covariant $\d-$modular forms (in
a sense generalizing that in \cite{difmod}). A natural question is
whether all $\d-$eigenforms can be obtained in this way. As we
shall see in this paper the answer is {\it no}.
 Indeed, we  provide, in this paper, a construction of
 $\d-$eigenforms $f^{\sharp}$ of weight $0$ and order $2$ that
 have  ``essentially the
same'' Hecke eigenvalues as certain classical  eigenforms $f$ of
weight $2$ (and order $0$). As we shall see, forms of weight $0$
(such as $f^{\sharp}$) are never $I-$linear combinations of forms
$f^{\phi^j}$. Having constructed the forms $f^{\sharp}$ one can
ask, of course, if any $\d-$eigenform is an $I-$linear combination
of forms $f^{\phi^j}$ and $(f^{\sharp})^{\phi^i}$; at this point
it is not clear what to expect.

The $\d-$Fourier expansion of $f^{\sharp}$ will be related in an
interesting way to the Fourier expansion of $f$. Indeed, if
$f=\sum a_n q^n$
 is a (classical) newform of weight
$2$ on $\Gamma_0(N)$
(which is not of ``CM type") with Fourier coefficients
$a_n \in \bZ$, then the $\d-$Fourier expansion of $f^{\sharp}$
will be a series $f^{\sharp}_{\infty}(q,q',q'')$ in $3$ variables
$q,q',q''$ which, after the substitution $q'=q''=0$, becomes
equal to the series
\[f^{(-1)}(q):=
\sum_{(n,p)=1} \frac{a_n}{n} q^n.\] (A similar, but more
complicated statement holds for $f$ of  ``CM type''.) The series
$f^{(-1)}$ is, of course, not the Fourier series of any
(classical) eigenform
 but, rather, a $p-$adic modular form in the sense of Serre; cf.
\cite{serre}, p. 115.  Note that, viewed as a function of elliptic
curves in the sense of Katz \cite{Katz} the $p-$adic modular form
$f^{(-1)}$ does not extend across the ``supersingular disks''
because, if this were  the case, $f^{(-1)}$ would define a
non-constant function on a projective modular curve. On the other
hand, remarkably, the $\d-$modular form $f^{\sharp}$ does extend
across the ``supersingular disks'' (this being the case with {\it
any} $\d-$modular form). One may ask if, in spite of this
phenomenon, $f^{\sharp}$ is, neverthless, a linear combination,
with isogeny covariant coefficients defined outside the
supersingular disks, of $\phi-$powers of $f^{(-1)}$; we will show
that this is not the case.

The idea in our construction of the forms $f^{\sharp}$ is to use
the Eichler-Shimura construction for the $f$'s in conjunction with
our theory of $\d-${\it characters} introduced in \cite{char}.
(Roughly speaking $\d-$characters are homomorphisms from the group
of $R-$points of an Abelian variety to the additive group of $R$
which, in coordinates, are given by expressions involving the
coordinates of the points and their iterated Fermat quotients.
They are arithmetic analogues of the {\it Manin maps}
 introduced by Manin in the context of the Mordell
conjecture over function fields
 \cite{man}.) Then our forms $f^{\sharp}$ will arise by composing
certain $\d-$characters of the modular Jacobians $J_1(N)$ with the
Abel-Jacobi maps $X_1(N) \ra J_1(N)$ that send  a fixed cusp into
$0$.

Here is the plan of this paper. In Section 2 we review (and
slightly extend) the concept of $\d-$modular form and Hecke
operators in \cite{difmod}.  Then we state one of our main results
about the existence  of the forms $f^{\sharp}$ and their
independence from $f$. In Section 3 we review results of
Eichler-Shimura and Manin-Drinfeld. Section 4 reviews
$\d-$characters \cite{char} and examines the existence of
eigenvectors in the space of $\d-$characters. In Section 5 we
conclude our construction of the forms $f^{\sharp}$; it will turn
out that the forms $f^{\sharp}$ ``vanish at all the cusps''. In
Section 6 we introduce $\d-${\it Fourier expansions} at $\infty$
and we compute them for our forms $f^{\sharp}$. In Section 7 we
use $\d-$Fourier expansions to prove, in particular,
 the independence of
$f^{\sharp}$ from $f^{(-1)}$. In Section 8 we use
 $\d-$Fourier expansions to
compute the effect of $\d-$Serre operators (in the sense of
\cite{book}) on $f^{\sharp}$. In Section 9 we prove that the
$\d-$Fourier expansion of $f^{\sharp}$ is in the domain of
definition of the (partially defined) {\it Hecke operator}
$T(p)_{\infty}$ and is an eigenvector of this operator. We  end
the paper by stating a result (whose proof will be given in a
subsequent paper \cite{dcc}) saying that $\d-$modular forms which
``vanish at the cusps'' and are in the domain of definition of
$T(p)_{\infty}$ automatically arise from composing $\d-$characters
of the modular Jacobians with Abel-Jacobi maps. This is, in some
sense, a converse of our existence results for the forms
$f^{\sharp}$ in the present paper.

{\bf Acknowledgement.} While writing this paper the author was
partially supported by NSF grant DMS 0552314.

\section{Main concepts}

\subsection{Prolongation sequences} Our main reference here
is  \cite{difmod}. We fix, throughout this paper, a prime integer
$p \geq 5$. Let $C_p(X,Y) \in \bZ[X,Y]$ be the polynomial with
integer coefficients \[C_p(X,Y):=\frac{X^p+Y^p-(X+Y)^p}{p}.\] A
$p-${\it derivation} from a ring $A$ into an $A-$algebra
$\varphi:A \ra B$ is a map $\d:A \ra B$ such that $\d(1)=0$ and
\[\begin{array}{rcl}
\d(x+y) & = &  \d x + \d y
+C_p(\varphi(x),\varphi(y))\\
\d(xy) & = & \varphi(x)^p \cdot \d y +\varphi(y)^p \cdot \d x
+p \cdot \d x \cdot \d y,
\end{array}\] for all $x,y \in A$. Given a
$p-$derivation we always denote by $\phi:A \ra B$ the map
$\phi(x)=\varphi(x)^p+p \d x$; then $\phi$ is a ring homomorphism.
A {\it prolongation sequence} is a sequence $S^*$ of  rings $S^n$,
$n \geq 0$, together with ring homomorphisms $\varphi_n:S^n \ra
S^{n+1}$ and $p-$derivations $\d_n:S^n \ra S^{n+1}$ such that
$\d_{n+1} \circ \varphi_n=\varphi_{n+1} \circ \d_n$ for all $n$.
We usually denote all $\varphi_n$ by $\varphi$ and all $\d_n$ by
$\d$ and we view $S^{n+1}$ as an $S^n-$algebra via $\varphi$. A
morphism of prolongation sequences, $u^*:S^* \ra \tilde{S}^*$ is a
sequence $u^n:S^n \ra \tilde{S}^n$ of ring homomorphisms such that
$\delta \circ u^n=u^{n+1} \circ \d$ and $\varphi \circ u^n=u^{n+1}
\circ \varphi$. Let $W$ be the ring of polynomials $\bZ[\phi]$ in
the indeterminate $\phi$. Then, for $w=\sum_{i=0}^r a_i \phi^i \in
W$, we set $deg(w):=\sum a_i$. If $a_r \neq 0$ we set $ord(w)=r$;
we also set $ord(0)=0$. For $w$ as above (respectively for $w \in
W_+:=\{\sum b_i \phi^i\ |\ \ b_i \geq 0\}$), $S^*$ a prolongation
sequence, and $x \in (S^0)^{\times}$ (respectively $x \in S^0$) we
can consider the element $x^w:=\prod_{i=0}^r \varphi^{r-i}
\phi^i(a)^{a_i} \in (S^r)^{\times}$ (respectively $x^w \in S^r$).
We let $W(r):=\{w \in W\ |\ ord(w) \leq r\}$.

Let $R:=R_p:=\hat{\bZ}_p^{ur}$ be the completion of the maximum
unramified extension of $\bZ_p$. Then $R$ has a unique
$p-$derivation $\d:R \ra R$ given by $\d x=(\phi(x)-x^p)/p$ where
$\phi:R \ra R$ is the unique lift of the $p-$power Frobenius map
on $R/pR$. One can consider the  prolongation sequence $R^*$ where
$R^n=R$ for all $n$. By a {\it prolongation sequence over $R$} we
understand a prolongation sequence $S^*$ equipped with a morphism
$R^* \ra S^*$. From now on all our prolongation sequences are
assumed to be over $R$.

\subsection{$\d-$modular forms} Our
main reference here is, again, \cite{difmod}. We fix, throughout
this paper, an integer $N \geq 1$, not divisible by $p$. For any
ring $S$  let us denote by $\bM(\Gamma_1(N),S)$ the set of all
triples $(E/S,\alpha,\omega)$ where $E/S$ is an elliptic curve,
$\omega$ is an invertible $1-$form on $E$, and
$\alpha:(\bZ/N\bZ)_S \ra E$ is a closed immersion of group schemes
(referred to as a $\Gamma_1(N)-${\it level structure}). Fix $w \in
W$ with $ord(w) \leq r$. A $\d-${\it modular form} of weight $w
\in W$ and order $r$ on $\Gamma_1(N)$ is  a rule $f$
  that associates to any prolongation sequence $S^*$
  of Noetherian, $p-$adically complete
  rings
  and any triple $(E/S^0, \alpha, \o) \in \bM(\Gamma_1(N),S^0)$
  an element $f(E/S^0, \alpha, \o, S^*) \in S^r$ such that the following
  properties are satisfied:

\begin{enumerate}

 \item $f(E/S^0, \alpha, \o, S^*)$  depends on the isomorphism class
  of  $(E/S^0, \alpha, \o)$ only.

\item  Formation of  $f(E/S^0, \alpha, \o, S^*)$  commutes with base change
$u^*:S^* \ra \tilde{S}^*$  i.e.
\[f(E \otimes_{S^0}
 \tilde{S}^0/\tilde{S}^0, \alpha \otimes \tilde{S}^0,
 u^{0*} \o, \tilde{S}^*)=
 u^r(f(E/S^0, \alpha, \o, S^*)).\]

\item $f(E/S^0, \alpha, \lambda \o, S^*)=\lambda^{-w} \cdot
 f(E/S^0, \alpha, \o, S^*)$
   for all $\lambda \in (S^0)^{\times}$.

\end{enumerate}

We denote by $M^r(\Gamma_1(N),R,w)$ the $R-$module of all
$\d-$modular forms over $R$ of weight $w \in W$ and order $r$ on
$\Gamma_1(N)$. Then the direct sum
\[M^r(\Gamma_1(N),R,*):=
\bigoplus_{w \in W(r)} M^r(\Gamma_1(N),R,w)\] has a natural
structure of graded ring. We view $M^{r}(\Gamma_1(N),R,*)$ as a
subring of $M^{r+1}(\Gamma_1(N),R,*)$ via $\varphi$ and we have
naturally induced ring homomorphisms $\phi:M^{r}(\Gamma_1(N),R,*)
\ra M^{r+1}(\Gamma_1(N),R,*)$ sending any $f \in
M^r(\Gamma_1(N),R,w)$ into $f^{\phi}:=\phi \circ f \in
M^{r+1}(\Gamma_1(N),R,\phi w)$; for $w=\sum a_i \phi^i \in W_+$ we
write $f^w:=\prod (f^{\phi^i})^{a_i}$. The rings
$M^r(\Gamma_1(N),R,*)$ are integral domains. Their union will be
denoted by $M^{\infty}(\Gamma_1(N),R,*)$.

We end our discussion here by noting that, by \cite{book}, p.252,
the spaces $M^r(\Gamma_1(N),R,w)$ embed into spaces
of {\it ordinary $\d-$modular forms}, denoted by
\[
M^r_{ord}(\Gamma_1(N),R,w)\] and defined exactly as the spaces
$M^r(\Gamma_1(N),R,w)$ except that instead of the set
$\bM(\Gamma_1(N),S^0)$ one considers the set
$\bM_{ord}(\Gamma_1(N),S^0)$ of all tuples in
$\bM(\Gamma_1(N),S^0)$ with ordinary reduction.

\subsection{$\d-$Hecke operators}
Again, our main reference here is \cite{difmod}. Assume $S^*$ is a
prolongation sequence of Noetherian, $p-$adically complete rings,
and let $\tilde{S}$ be a finite \'{e}tale over-ring of $S^0$.
 Then, by \cite{difmod}, (1.6), there is a unique structure
of prolongation sequence on $S^*\otimes_{S^0} \tilde{S}:= (S^n
\otimes_{S^0} \tilde{S})$ compatible (in the obvious sense) with
that of $S^*$. Now let $l$ be a prime integer not dividing $Np$.
 Let $f \in M^r(\Gamma_1(N),R,w)$ be a
$\d-$modular form. We can define a $\d-$modular form $T(l) f \in
M^r(\Gamma_1(N),R,w)$ by the formula
\begin{equation}
\label{tofl} (T(l) f)(E/S^0, \alpha, \o, S^*)= \sum_{i=0}^l
f(\tilde{E}_i/\tilde{S}, u_i \circ \alpha, u_{i*} \o, S^* \otimes
_{S^0} \tilde{S}) \end{equation}
 where $\tilde{S}$ is any finite
\'{e}tale over-ring of $S^0$ such that the group scheme of points
of order $l$ of $\tilde{E}:=E \otimes_{S^0} \tilde{S}$
 is isomorphic to  $({\bf Z}/l \bZ)^2_{\tilde{S}}$
(hence the elliptic curve $\tilde{E}$ has exactly $l+1$ finite,
flat subgroup schemes $H_0,...,H_l$ of rank $l$),
$\tilde{E}_i=\tilde{E}/H_i$, $u_i:\tilde{E} \ra \tilde{E}_i$ are
the natural projections, and the $u_{i*} \o$'s are induced by $\o$
via pull-back to $\tilde{E}$ followed by trace to the
$\tilde{E}_i$'s. In the above we can always assume $\tilde{S}$ is
Galois over $S^0$. Note that
 $(T(l) f)(E/S^0,\alpha, \o, S^*)$ which is,
 a priori, an element of $S^r \otimes_{S^0} \tilde{S}$,
 actually belongs to $S^r$, and does not depend on the choice
 of $\tilde{S}$.

We refer to the maps $T(l):M^r(\Gamma_1(N),R,w) \ra
M^r(\Gamma_1(N),R,w)$ as {\it $\d-$Hecke operators}. Clearly these
maps commute with $\phi$.
 For $r=0$ and $w=m \in \bZ$ one can normalize our $T(l)$
in the classical fashion by
 considering the operators
 \[T_m(l):=l^{m-1}T(l):M^0(\Gamma_1(N),R,m) \ra
 M^0(\Gamma_1(N),R,m).\]

\subsection{Classical eigenforms} Our main references here are
 \cite{RS, con}.
Denote
 by $S_m(\Gamma_1(N),\bC)$  the space of (classical)
  cusp forms of weight $m$ on $\Gamma_1(N)$ over the complex field $\bC$.
   On this space one has
Hecke operators $T_m(n)$ acting, $n \geq 1$. An {\it eigenform} $f
\in S_m(\Gamma_1(N),\bC)$ is a nonzero element which is a
simultaneous eigenvector for all $T_m(n)$, $n \geq 1$. An
eigenform $f=\sum_{n \geq 1} a_n q^n$, $a_n=a_n(f)$, is {\it
normalized} if $a_1=1$; in this case $T_m(n)f=a_n \cdot f$ for all
$n \geq 1$. One associates to any  eigenform $f \in
S_m(\Gamma_1(N),\bC)$ its {\it system of eigenvalues} $l \mapsto
a_l$, $(l,N)=1$. A {\it newform} is a normalized eigenform whose
 system of eigenvalues does not come from
 a system of eigenvalues associated
to an eigenform in $S_m(\Gamma_1(M),\bC)$ with $M\ |\ N$, $M \neq
N$. For any normalized eigenform $f \in S_m(\Gamma_1(N),\bC)$ one
may consider the subring $\cO_f$ of $\bC$ generated by all
$a_n(f)$, $n \geq 1$; then $\cO_f$ is a finite $\bZ-$algebra  and
one denotes by $K_f$ its fraction field. If $Q \geq 1$ is any
integer we denote by $\cO_f^{(Q)}$ the subring of $\bC$ generated
by all $a_l(f)$, where $l$ is prime, not dividing $Q$.

We will later need to consider the subspace $S_m(\Gamma_0(N),\bC)$
of $S_m(\Gamma_1(N),\bC)$ of all cusp forms of weight $m$ on
$\Gamma_0(N)$. Recall that if $f=\sum a_n q^n$ is an eigenform in
$S_m(\Gamma_0(N),\bC)$ then
\begin{equation} \label{tzu}
\begin{array}{rcl}
a_{n_1n_2} & = & a_{n_1}a_{n_2}\ \
for\ \  (n_1,n_2)=1\\
a_{l^{i-1}}a_l & = & a_{l^i}+l^{m-1}
a_{l^{i-2}}\ \ for \ \ l\ \  prime,\ \
(l,N)=1\ \ and\ \ i \geq 2,\\
a_{l^{i-1}}a_l & = & a_{l^i}\ \ for \ \ l \ \ prime,\ \
l|N\ \ and\ \ i \geq 2.\end{array} \end{equation}

\subsection{$\d-$eigenforms}
A non-zero $\d-$modular form  $h \in M^r(\Gamma_1(N),R_p,w)$ is
called a $\d-${\it eigenform} if $T(l)h=\lambda_l \cdot h$,
$\lambda_l \in R_p$, for all primes $l$ not dividing $Np$. A
$\d-$eigenform is said to {\it belong} (outside $Np$) to a
(classical) normalized eigenform $f=\sum a_n q^n \in
S_m(\Gamma_1(N),\bC)$ if there exist a (necessarily injective)
ring homomorphism $\e:\cO_f^{(Np)} \ra R_p$ and an integer $e \in
\bZ$ such that
$\lambda_l=l^e  \e(a_l)$
 for all primes
$l$ not dividing $Np$. We then say that $f^{\sharp}$ belongs to
$f$ with {\it character} $\e$ and {\it exponent} $e$.

Note that $\e$ is uniquely determined by $f^{\sharp}$. Indeed
assume $\e,\e':\cO_f^{(Np)} \ra R_p$ are ring homomorphisms and
$e,e'$ are integers such that
$l^e  \e(a_l)=l^{e'}
\e'(a_l)$ for all primes $l$ not dividing $Np$ and assume $\e
\neq \e'$. Then clearly $e \neq e'$. Set $L:=l^{e-e'} \neq 1$ and,
since $\e \neq \e'$, one can choose a prime $l$ not dividing $Np$
such that $a_l \neq 0$. Let $\Phi(t)=t^d+b_1^{d-1}+...+b_d \in
\bQ[t]$ be the minimal polynomial of $a_l$ over $\bQ$. Then both
$\e(a_l)$ and $L \cdot \e(a_l)$ are roots of $\Phi(t)$. Hence
$\e(a_l)$ is a root of
$\Psi(t):=\Phi(t)-L^{-d}\Phi(Lt)=\sum_{i=1}^d
b_i(1-L^{-i})t^{d-i}$. Since $\Psi(t)$ has degree $\leq d-1$ we
must have $\Psi(t)=0$ hence $\Phi(t)=t$ hence $a_l=0$, a
contradiction.

\subsection{$\d-$eigenforms arising from classical eigenforms}
There is an ``easy"  way to construct $\d-$ eigenforms belonging
to classical eigenforms $f$ by taking linear combinations of
``$\phi-$ powers of $f$" with {\it isogeny covariant} $\d-$modular
forms (in a sense slightly generalizing that in \cite{difmod}). In
what follows we explain this construction. We should point
 out that the forms $f^{\sharp}$ mentioned in the Introduction
will be shown not to be obtainable via this construction.

Let $F \in M^r(\Gamma_1(N),R,w)$ be a $\d-$ modular form of weight
$w=\sum n_i \phi^i$ on $\Gamma_1(N)$. Assume $deg(w):=\sum n_i$ is
even. Generalizing the level one definition in \cite{difmod} we
say that $F$ is {\it isogeny covariant} if  for any prolongation
sequence $S^*$, any triples
$(E_1,\alpha_1,\omega_1),(E_2,\alpha_2,\omega_2) \in
\bM(\Gamma_1(N),S^0)$, and
 any isogeny $u:E_1 \ra E_2$  of degree prime to $p$,
 with $\omega_1=u^* \omega_2$ and $u \circ \alpha_1=\alpha_2$, we have
\[F(E_1,\alpha_1,\omega_1,S^*)=deg(u)^{-deg(w)/2} \cdot F(E_2,\alpha_2,
\omega_2,S^*).\]

\begin{example}
By  \cite{book}, p. 268 and Theorem 8.83,
 for each $r \geq 1$  the $R_p-$module of
isogeny covariant $\d-$modular forms in
$M^r(\Gamma_1(N),R_p,-1-\phi^r)$ is free of rank one. Following
\cite{book} we shall denote by $f^r=f^r_{crys}$ a basis of this
rank one module. (So the upper $r$ is an index, not an exponent.
Recall from \cite{book} that $f^r$ is constructed via crystalline
cohomology.)
\end{example}

We denote by $\cI \subset M^{\infty}(\Gamma_1(N),R,*)$ the
multiplicative system of all non-zero isogeny covariant
$\d-$modular forms and by $\cJ \subset \cI$ the multiplicative
system generated by all $(f^r)^{\phi^s}$ for $r \geq 1$ and $s
\geq 0$. The $R-$linear spans of $\cI$ and $\cJ$ will be denoted
by $I$ and $J$ respectively. Then $I$ is a ring, $J$ is a subring
 of $I$, and it is tempting to conjecture \cite{Barcau, siegel}
that $ J \otimes \bQ=I \otimes \bQ. $

\begin{lemma}
\label{ploua} If $F \in M^r(\Gamma_1(N),R,w)$ is isogeny covariant
and \[G \in  M^r(\Gamma_1(N),R,v)\] is any $\d-$modular form then,
for any prime $l$ not dividing $Np$,
\[T(l)(F \cdot G)=l^{-deg(w)/2} \cdot F \cdot T(l)G.\]
In particular, if $G$ is a $\d-$eigenform belonging to the
classical normalized eigenform $f \in S_m(\Gamma_1(N),\bC)$ with
character $\e$ and exponent $e$ then $F \cdot G$ is a
$\d-$eigenform belonging to $f$ with character $\e$ and exponent
$e-\frac{deg(w)}{2}$.
\end{lemma}

\begin{proof} This follows from
a computation similar to the one in \cite{char}, p.125 (where the
case $N=1$ was treated).
\end{proof}

 Now let $f \in
S_m(\Gamma_1(N),\bC)$, $f=\sum a_n q^n$,  be a normalized
eigenform of weight $m$ and  let $\rho:\cO_f [1/N,\zeta_N] \ra
R_p$ be any ring homomorphism, where $p$ does not divide $N$.
Then,
 by the ``$q-$expansion principle'' \cite{DI}, pp. 70 and 112,
 $f$ naturally defines (via $\rho$) a rule $f^{\rho}$
(compatible with base change
 and homogeneous of degree $-m$)
 that attaches to any $R_p-$algebra
 $S$ and any triple $(E/S,\alpha,\omega) \in \bM(\Gamma_1(N),S)$ an
 element $f^{\rho}(E/S,\alpha,\omega) \in S$ depending only on the
 isomorphism class of the triple. (Here it is essential that we have a fixed
 primitive $N-$th root of unity $\rho(\zeta_N)$ in $R_p$.)
 Then $f^{\rho}$ induces
   a $\d-$modular form (still denoted by) $f^{\rho} \in
M^0(\Gamma_1(N),R_p,m)$, of weight $m$, defined by the formula
\[f^{\rho}(E/S^0,\alpha,\omega,S^*):=f^{\rho}(E/S^0,\alpha,\omega).\]
  The
composition $f^{\rho \phi^j}:=\phi^j \circ f^{\rho}$ is a well
defined element of $M^j(\Gamma_1(N),R_p,\phi^jw)$.  There is an
obvious compatibility between the classical and our Hecke
operators $T_m(l)$, which yields:
\begin{equation}
\label{jeau}
\begin{array}{rcl}
 T(l) f^{\rho \phi^j} &
= & (T(l)f^{\rho})^{\phi^j}\\
\  & = & (l^{1-m}T_m(l)f^{\rho})^{\phi^j}\\
\  & = & l^{1-m} (T_m(l)f)^{\rho\phi^j}\\
\  &  = & l^{1-m} \phi^j (\rho(a_l)) \cdot f^{\rho \phi^j}
\end{array} \end{equation}
 for all primes $l$ not
dividing $Np$. So we see that $f^{\rho \phi^j}$ is a
$\d-$eigenform of order $j$ and  weight $m\phi^j$ which belongs to
$f$ with character $\phi^j \circ \rho$ and exponent $e=1-m$.

In particular, according to Lemma \ref{ploua},
if $\phi^b \circ \rho=\rho$ for some $b \geq 1$, and if $a \geq 0$, then
any $\d-$modular function of weight $w$ in the $I-$linear span of
\[\{
f^{\rho\phi^a},f^{\rho\phi^{a+b}},f^{\rho\phi^{a+2b}},f^{\rho\phi^{a+3b}},
...\}\]
is a $\d-$eigenform belonging to $f$ with character $\phi^a \circ \rho$
and exponent
 \begin{equation}
 \label{bum}
e=1-\frac{m+deg(w)}{2}.
 \end{equation}

\begin{lemma}
\label{langgg}
 Let $f \in S_2(\Gamma_1(N),\bC)$ be a
normalized eigenform and let $\tilde{f}$
be a non-zero
$\d-$modular form of weight $0$. Then $\tilde{f}$ cannot be
 in the
$I-$linear span of
\[\{f^{\rho},f^{\rho \phi},f^{\rho\phi^2},f^{\rho \phi^3},...\}.\]
\end{lemma}

\begin{proof}
Assume the conclusion is false. We may assume
$\tilde{f}=\sum F_a \cdot f^{\rho \phi^a}$,
where $F_a$ are isogeny covariant of weight $-2\phi^a$.
To get a contradiction we need to check the following:

{\it Claim.  For any $0 \leq a \leq r$ there are no non-zero
isogeny covariant $\d-$modular forms in
$M^r(\Gamma_1(N),R_p,-2\phi^a)$.}

 We fix $a$ and prove this
claim by induction on $r$. Assume first $r=a$. If $a=0$ then the
claim follows from \cite{book}, Proposition 8.75. If $a \geq 1$
our claim  follows from \cite{book}, Theorem 8.83, assertion 2. To
perform the induction step assume $r>a$. Then, by \cite{book},
Corollary 8.40 and  Proposition 8.75, if $h \in
M^r(\Gamma_1(N),R_p,-2\phi^a)$ is isogeny covariant then
$\partial_r h=0$ where $\partial_r$ is the $\d-$Serre operator in
loc. cit. This easily implies that $h \in
M^{r-1}(\Gamma_1(N),R_p,-2\phi^a)$ and we conclude by the
induction hypothesis. (The various results in \cite{book} quoted
above apply to our situation in view of \cite{book}, Proposition
8.22. Note also that the special case $N=1$ of our claim was
proved by Barcau \cite{Barcau}.)
\end{proof}

\subsection{Ordinary
$\d-$modular forms arising from $p-$adic modular forms}

The main refernces here are \cite{serre, Katz, goren}.
Let $g \in \bQ_p[[q]]$ be a $p-$adic modular form of weight $m \in \bZ$ in
the sense of Serre. Fix a homomorphism
$\rho:\bZ[\zeta_N,1/N] \ra R_p$. Then $g$ induces
a $p-$adic modular form $g^{\rho}$ of level $N$,
weight $m$ and growth $1$ in the sense of Katz \cite{Katz};
cf. \cite{goren}, Theorem 6.21, p. 158.
On the other hand $g^{\rho}$ induces an ordinary
$\d-$modular form (still denoted by)
$g^{\rho} \in M^0_{ord}(\Gamma_1(N),
R_p,m)$. So for each $j \geq 0$ we may consider the ordinary
$\d-$modular form $g^{\rho \phi^j} \in M^j_{ord}(\Gamma_1(N),R_p,m\phi^j)$.

In particular, if $f=\sum a_n q^n \in S_2(\Gamma_1(N),\bC)$
is a normalized eigenform of weight $2$ with $a_n \in \bZ$ then,
by \cite{serre}, p. 115, the series
\[f^{(-1)}=\sum_{(n,p)=1}\frac{a_n}{n} q^n\]
is a $p-$adic modular form of weight $0$ such that
$T_0(l)f^{(-1)}=l^{-1}a_lf^{(-1)}$ for $l$ prime different from $p$.
It immediately follows that \[f^{(-1)\rho\phi^j} \in
M^j_{ord}(\Gamma_1(N),R_p,0)
\] is, in the obvious sense,
 an ordinary $\d-$eigenform
belonging to $f$ with exponent $0$. So any $R-$linear
combination of such forms will have the same property.

Note that, if in the definition of isogeny covariant $\d-$modular forms,
one replaces $\bM(\Gamma_1(N),S^0)$ by
$\bM_{ord}(\Gamma_1(N),S^0)$ one obtains the notion of {\it ordinary
isogeny covariant $\d-$modular form}. Let ${\mathcal I}_{ord}$
be the multiplicative system of all such forms
and let $I_{ord}$ be the $R-$linear span of ${\mathcal I}_{ord}$.
Then $I_{ord}$ is a ring.

\begin{lemma}
\label{mabatelakap}
Let $f \in S_2(\Gamma_1(N),\bC)$
be a normalized eigenform of weight $2$ with $a_n \in \bZ$
and let $\tilde{f}$ be an ordinary $\d-$modular form
of weight $0$ which is in the  $I_{ord}-$linear span of the set
\[
\{f^{(-1)\rho}, f^{(-1)\rho \phi}, f^{(-1)\rho \phi^2}, f^{(-1)\rho
\phi^2},...\}.
\]
Then $\tilde{f}$ is in the $R-$linear span of this set.
\end{lemma}

\begin{proof}
Let $\tilde{f}=\sum_j F_j \cdot f^{(-1) \rho \phi^j}$, $F_j \in I_{ord}$.
Picking the weight $0$ components we may assume $F_j$ have weight $0$.
So we are reduced to showing that any ordinary isogeny covariant
$\d-$modular form of weight $0$ is a constant in $R$. This follows
from \cite{cr97}, Propositions 7.21 (plus the Remark after it) and
7.23.
\end{proof}

\subsection{The forms $f^{\sharp}$}
The main purpose of this paper is to provide a construction  of
 $\d-$ eigenforms $f^{\sharp}$
 of weight $0$ and order $2$ belonging to classical
eigenforms $f$ of weight $2$. As we shall see,
the forms $f^{\sharp}$ will
be neither $I-$linear combinations
of $\phi-$powers of $f$ nor $I_{ord}-$linear combinations of
$\phi-$powers of $f^{(-1)}$.
Here is one of our main results.
 This result will be complemented by other results
later in the paper; cf. Remark \ref{ccci} below.

\begin{theorem}
\label{main} Let   $f \in S_2(\Gamma_1(N),\bC)$
 be a   newform of
weight $2$  on $\Gamma_1(N)$, $N>4$, and let $g:=[K_f:\bQ]$. Then,
for any sufficiently large prime $p$, and any embedding
$\rho:\cO_f[1/N,\zeta_N] \ra R_p$,
 there exist
$\d-$eigenforms
\[f_1^{\sharp},...,f^{\sharp}_g \in M^2(\Gamma_1(N),R_p,0),\]
of weight $0$ and  order $2$ on $\Gamma_1(N)$, such that

1)  $f^{\sharp}_j$ belongs to $f$ with exponent $0$,

2) $f^{\sharp}_j$ is not in the $I-$linear span of
$\{f^{\rho},f^{\rho \phi},f^{\rho\phi^2},f^{\rho \phi^3},...\}$,

3) $f^{\sharp}_1,...,f^{\sharp}_g$ are $R_p-$linearly independent.
\end{theorem}

\begin{remark}
\label{ccci}

1) Assertion 2
 follows directly from Lemma \ref{langgg}.

2) As we shall see, $f_j^{\sharp}$ themselves should be morally
viewed as a kind of ``$\d-$cusp forms'' in the sense that they
``vanish at the cusps''; cf. Remark \ref{ccuusspp}.

3) In case $g=1$ (i.e. $a_n(f) \in \bQ$, equivalently $a_n(f) \in
\bZ$) the form $f^{\sharp}=f_1^{\sharp}$ will be essentially
canonically associated to $f$ and its $\d-${\it Fourier expansion}
will be closely related to that of $f$; cf.
Theorems \ref{fur} and
\ref{fur2}.
We will show that $f^{\sharp}$ is not an $R_p-$linear combination
of $\phi-$powers of $f^{(-1)}$; cf. Theorem \ref{doii}.
 Also we shall compute the
effect of the $\d-${\it Serre operators}
on $f^{\sharp}$; cf. Propositions \ref{se1} and \ref{se2}.

4) It would be interesting to extend the above Theorem to $f$'s of
higher weight.

5) Let $g=1$, $f^{\sharp}:=f^{\sharp}_1$,
and  $0 \leq a < b \leq r$, $r \geq 2$;
then the $\d-$eigenforms
$(f^a)^{\phi^{b-a}} \cdot
f^{\sharp}$ and $f^a \cdot f^b \cdot f^{\rho}$ of
order $r$ have the same  weight $-\phi^a-\phi^b$ and belong to $f$
with the same character and same
exponent $e=1$. One can ask if these forms are
$R_p-$linearly dependent. The answer is negative as one will see
in Theorem \ref{doii}.
\end{remark}

\section{Review of Eichler-Shimura and Manin-Drinfeld}

We need to review some basic facts about modular curves.
 The
references for this section are \cite{KM, DI, RS, con}.
 Fix an integer $N \geq 4$. Then the {\it modular curve}
\begin{equation}
\label{Y1(N)} Y_1(N)/\bZ[1/N]
\end{equation}
 is the scheme whose $S-$points
($S$ any scheme over $\bZ[1/N]$) identify with isomorphism classes
of pairs $(E,\alpha)$ where $E/S$ is an elliptic curve and
$\alpha:(\bZ/N\bZ)_S \ra E$ is a closed immersion of group
schemes. Recall that $Y_1(N)$ is smooth affine of relative
dimension one over $\bZ[1/N]$, with geometrically irreducible
fibers.

Let now $l$ be a prime integer not dividing $N$.
 There is a scheme
\begin{equation}
\label{Y1(N,l)} Y_1(N,l)/\bZ[1/Nl]
\end{equation}
  whose $S-$points
identify with
triples $(E,\alpha,H)$ where $(E,\alpha)$ is as above and
$H$ is a finite flat subgroup scheme of $E$ of rank $l$.
 Recall that $Y_1(N,l)$ is smooth affine of relative dimension one
over $\bZ[1/Nl]$, with geometrically irreducible fibers. Over
$\bZ[1/Nl]$ one can consider the natural projections
\begin{equation}
\label{bfcal}
\begin{array}{rcl}
 \i, \s:Y_1(N,l) & \ra & Y_1(N)\\
\i(E,\alpha,H) & = & (E,\alpha),\\
\s(E,\alpha,H) & = & (E/H,u \circ \alpha),\end{array}
\end{equation} where $u:E \ra E/H$ is the
canonical projection. Moreover $\i$ and $\s$ are \'{e}tale above
$\bZ[1/Nl]$. For details on the discussion above see  \cite{KM}
pp. 87, 117, 125, 129, \cite{DI}, pp. 69-72, \cite{con}, pp.
212-213. (For $\sigma_1$ we use the convention in \cite{con}
rather than that in \cite{DI}.)

In what follows we need to consider the ``compactified" situation.
We assume $N>4$.  Recall that $Y_1(N)$ is an open set in its {\it
Deligne-Rapoport compactification}
$X_1(N)/\bZ[1/N]$
 which is a proper smooth
scheme.   Cf. \cite{DI}, p. 79. The complex points of
$X_1(N)
\backslash Y_1(N)$
 are called the {\it cusps} of $X_1(N)$ and they come from
 $\bZ[1/N,\zeta_N]-$points of
 $X_1(N)$. As
usual we denote by $J_1(N)$  the Jacobian of $X_1(N)$ viewed as an
Abelian scheme over $\bZ[1/N]$. Let $X_1(N,l)_{\bC}$ be a smooth
complete model over $\bC$ of the complex curve $Y_1(N,l)_{\bC}$;
the morphisms in Equation \ref{bfcal} induce morphisms
\[\sigma_1,\sigma_2:X_1(N,l)_{\bC} \ra X_1(N)_{\bC}.\]
These morphisms induce endomorphisms $T(l)_*$ of $J_1(N)_{\bC}$ as
follows: if $D$ is a divisor of degree $0$ on $X_1(N)_{\bC}$ and
$[D]$ is the point of $J_1(N)_{\bC}$ representing $D$ then
\begin{equation}
\label{forforT} T(l)_*[D]:=[\sigma_{2*}\sigma_1^* D].
\end{equation}
The endomorphisms $T(l)_*$  have models over $\bZ[1/N]$ (arising
from N\'{e}ron model theory); cf. \cite{con}. We will need the
following basic construction due to Eichler-Shimura (cf.
\cite{con}, p. 215, \cite{DS}, pp. 241-242):

\begin{theorem}
\label{eichshi} (Eichler-Shimura) Let $f \in
S_2(\Gamma_1(N)),\bC)$ be a newform.  Then there exists an Abelian
variety $A=A_{\bQ}$ defined over $\bQ$  of dimension $[K_f:\bQ]$,
a ring homomorphism $\iota: \cO_f \ra End(A/\bQ)$, and a dominant
homomorphism $\pi:J_1(N)_{\bQ} \ra A$ defined over $\bQ$ such that
the following hold:

1) For all primes $l$ we have a commutative diagram
\[\begin{array}{ccc}
J_1(N)_{\bQ} & \stackrel{T(l)_*}{\longrightarrow} & J_1(N)_{\bQ}\\
\pi \downarrow & \  & \downarrow \pi\\
A & \stackrel{\iota(a_l(f))}{\longrightarrow} & A
\end{array}
\]

2) The image of the pull-back map
\[\pi^*:H^0(A_{\bC},\Omega) \ra
H^0(J_1(N)_{\bC},\Omega) \simeq S_2(\Gamma_1(N))_{\bC}\] is the
$\bC-$linear span of $\{f^{\sigma}\ |\ \sigma:K_f \ra \bC\}$.

3) If $f \in S_2(\Gamma_0(N),\bC)$, $g=1$,  and $p$ is a
sufficiently big prime then $a_p$
 equals the trace of the $p-$power Frobenius on the elliptic
 curve obtained by reducing $A$ mod $p$.
\end{theorem}

\begin{remark}
\label{CM} If $g=1$ then $A$ in the above Theorem is an elliptic
curve over $\bQ$ (which, by condition 2 in the Theorem,  is
uniquely determined by $f$ up to isogeny); we say that $f$ is of
CM type or not of CM type according as $A$ has CM or does not have
CM. Assertion 3 in the above Theorem, appropriately reformulated,
holds for $\Gamma_1(N)$ and arbitrary $g$; for our purposes here
we will not need this more general statement.
\end{remark}

 On the other hand we will need the following theorem due
to Manin and Drinfeld (cf. \cite{man1} or \cite{lang}, p. 62):

\begin{theorem}
\label{mandrin} (Manin-Drinfeld) If $D$ is a divisor of degree $0$
on $X_1(N)_{\bC}$ supported in the set of cusps then its image
$[D]$ in $J_1(N)_{\bC}$ is a torsion point.
\end{theorem}

\section{$\d-$characters}

We start by reviewing some concepts from  \cite{difmod,
char}. A $\d-${\it morphism} of order $r$, $f:X \ra Y$,
between two $R-$schemes is a rule that attaches to any
prolongation sequence $S^*$ of $p-$adically complete rings
 a map $f_{S^*}:X(S^0) \ra
Y(S^r)$ which is functorial in $S^*$.
 In the special case when
 $X$ is smooth over $R$ and $Y=\bA^1$ is the
affine line any $\d-$morphism $f:X \ra \bA^1=\bA^1_R$ is
completely determined by the map $f_{R^*}:X(R) \ra \bA^1(R)=R$. We
denote by $\cO^r(X)$ the ring of all $\d-$morphisms $X \ra \bA^1$
of order $r$. Assume
 that $G$ is a smooth group scheme over $R$. A $\d-$morphism
$f:G \ra \bA^1={\bf G}_a$ for  which $f_{R^*}:G(R) \ra {\bf
G}_a(R)=R$ is a group homomorphism into the additive group of $R$
is called a $\d-${\it character}.
We denote by $\bX^r(G)$ the $R-$module of $\d-$characters of $G$
of order $r$.

\begin{theorem}
\label{mythm} \cite{char} Let $A$ be an Abelian scheme over $R$ of
relative dimension $g$. Then
$(r-1)g \leq rank_R \bX^r(A) \leq  rg$.
\end{theorem}

\begin{proof}
This is contained in  \cite{char}, pp. 325-326.
\end{proof}

As explained in \cite{char}, the $\d-$characters $\psi:A \ra {\bf G}_a$
of an Abelian variety should be viewed as
arithmetic analogues of the Manin maps in \cite{man}; Manin maps
are homomorphisms $A(L) \ra L$ defined for any Abelian variety $A$
over
 a field $L$ of characteristic zero equipped
with a non-zero derivation $D:L \ra L$. In our theory the role of
the derivation $D$ is
 played by the $p-$derivation $\d$.

\begin{remark}
For $N \geq 4$ the description of points of $Y_1(N)$ immediately
implies that we have an identification
\begin{equation}
\label{identif1} M^r(\Gamma_1(N),R,0) \simeq \cO^r(Y_1(N)_R).
\end{equation}
  More generally, the spaces $M^r(\Gamma_1(N),R,w)$ and
$M^r_{ord}(\Gamma_1(N),R,w)$ identify with the spaces
$M^r_{Y_1(N)_R}(w)$ and $M^r_{Y_1(N)_{R,ord}}(w)$ in \cite{book},
p. 251, respectively.
\end{remark}

The above suggests the following:

\begin{definition}
\label{cussp} A $\d-$modular form $f \in M^r(\Gamma_1(N),R,0)
\simeq \cO^r(Y_1(N)_R)$ is called {\it $\d-$holomorphic} if it
lies in the image of
$\cO^r(X_1(N)_R) \ra \cO^r(Y_1(N)_R)$.
If  $\rho:\bZ[1/N,\zeta_N] \ra R$ is an embedding then the {\it
cusps} of $X_1(N)_R$ (with respect to $\rho$) are the $R-$points
of $X_1(N)_R$ obtained as  images, via $\rho$, of the cusps of
$X_1(N)_{\bC}$ (viewed as $\bZ[1/N,\zeta_N]-$points). A
$\d-$modular form $f \in M^r(\Gamma_1(N),R,0)$ is a {\it $\d-$cusp
form} (with respect to $\rho$) if it is $\d-$holomorphic and it
vanished at all the cusps of $X_1(N)_R$ (with respect to $\rho$).
\end{definition}

\begin{remark}
\label{tutu} For the next Proposition and its proof it is useful
to introduce some  terminology and record some facts about it. If
$K$ is a field, $V$ is an $n-$dimensional $K-$linear space, and
$T \in End(V)$ then, by the {\it eigenvalues} of $T$ on $V$ we
mean the eigenvalues (in an algebraic closure $K^a$ of $K$) of any
matrix in $Mat_n(K)$ representing $T$. If $W \subset V$ is a
subspace with $TW \subset W$ then all eigenvalues of $T$ on $W$
and  all eigenvalues of $T$ on $V/W$ are also  eigenvalues of $T$
on $V$. If $L$ is a field extension of $K$ we say that $T \in
End(V)$ is {\it diagonalizable} over $L$ if $V \otimes_K L$ has an
$L-$basis consisting of eigenvectors of $T$. More generally,
$T_1,...,T_s \in End(V)$ are said to be {\it simultaneously
diagonalizable} over $L$ if $V \otimes_K L$ has an $L-$basis
consisting of common eigenvectors of $T_1,...,T_g$. We will use
the following trivial facts:

1) $T$ is diagonalizable over $K^a$ if and only if the minimal
polynomial of $T$ on $V \otimes_K K^a$  has simple roots only.

2) If $T$ is diagonalizable over $K^a$ and all its eigenvalues are
in $K$ then $T$ is diagonalizable over $K$.

3) If $T_1,...,T_s$ are pairwise commuting and each of them is
diagonalizable over $K$ then $T_1,...,T_s$ are simultaneously
diagonalizable over $K$.
\end{remark}

\begin{proposition}
\label{eigch} Let $A/R$ be an Abelian scheme of relative dimension
$g$ and $\tau_1,...,\tau_s \in End(A/R)$ be commuting
endomorphisms each of which annihilates a polynomial with
$\bZ-$coefficients with simple complex roots only. Let ${\mathcal
T}$ be the subring of $End(A/R)$ generated by $\tau_1,...,\tau_s$
and let $K$ be the fraction field of $R$. Assume the tangent maps
$d \tau_1,...,d \tau_s \in End(Lie(A_K/K))$ have all their
eigenvalues in $K$. (Here we view $Lie(A_K/K)$ as the tangent
space to $A_K$ at the origin.) Then there exist $R-$linearly
independent  $\d-$characters $\psi_1,...,\psi_g:A \ra {\bf G}_a$
of order $2$ and  ring homomorphisms $\e^1,...,\e^g:{\mathcal T}
\ra R$ such that $\psi_j \circ \tau=\e_j(\tau) \cdot \psi_j$ for
all $\tau \in {\mathcal T}$, $j=1,...,g$.
\end{proposition}

\begin{proof}
Let the formal group $A^{for}$ of $A$  be identified with  $Spf\
R[[x]]$ where $x=\{x_1,...,x_g\}$ are some variables. Let
$L=(L_1,...,L_g) \in K[[x]]^g$ be the logarithm of the formal
group law of $A$ with respect to $x$. Let $x'$ and $x''$ be
additional $g-$tuples of variables and let
$K[[x]] \stackrel{\phi}{\ra} K[[x,x']] \stackrel{\phi}{\ra}
K[[x,x',x'']]$ be the ring homomorphisms, extending
$\phi:R \ra R$, defined by
$\phi(x)=x^p+px'$, $\phi(x')=(x')^p+px''$. By Theorem \ref{mythm}
we have $rank_R\bX^2(A) \geq g$. Note that $End(A/R)^{op}$ acts on
$\bX^2(A)$ (and hence on $\bX^2(A) \otimes K$) by the formula
$(\tau,\psi) \mapsto \tau^* \psi:=\psi \circ \tau$, $\tau \in
End(A/R)$. On the other hand $End(A/R)^{op}$ naturally acts on
$R[[x]]$ via $(\tau,x) \mapsto \tau^* x$ and one can extend this
action uniquely to an action on
$K[[x,x',x'']]=K[[x,\phi(x),\phi^2(x)]]$ such that the induced
endomorphisms $\tau^*$ on the latter ring satisfy
$\phi(\tau^* G)=\tau^*(\phi(G))$
for all $G \in K[[x,x']]$.
 By
\cite{char}, Lemma 2.8, there exists an injective $K-$linear
$End(A/R)^{op}-$equivariant map
$\bX^2(A) \otimes K \ra K[[x,x',x'']]$
whose image lies in the $K-$linear space
\begin{equation}
\label{spnn} V:=\sum_{i=1}^g(K \cdot L_i+K \cdot \phi(L_i)+K \cdot
\phi^2(L_i)). \end{equation}
 The logarithm $L:A_K^{for} \ra ({\bf
G}_{a,K}^{for})^g$ is an isomorphism of formal groups over $K$
hence, in particular, the endomorphisms $\tau_i:A^{for} \ra
A^{for}$ induce, via $L$, endomorphisms of $({\bf
G}_{a,K}^{for})^g$, hence matrices $M(\tau_i)=(m_{js}(\tau_i)) \in
Mat_g(K)$ satisfying
\begin{equation}
\label{8347t} \tau^*_i L_j=\sum_{s=1}^g m_{js}(\tau_i) L_s.
\end{equation}
On the other hand, taking the $Lie$ functor we get commutative
diagrams
\[\begin{array}{ccccc}
Lie(A_K/K) & \simeq & Lie(A_K^{for}) & \stackrel{dL}{\ra} &
Lie({\bf G}_{a,K}^{for})^g\\
d \tau_i \downarrow & \  & d \tau_i \downarrow & \  & \downarrow
M(\tau_i)\\
Lie(A_K/K) & \simeq & Lie(A_K^{for}) & \stackrel{dL}{\ra} &
Lie({\bf G}_{a,K}^{for})^g
\end{array}
\]
showing that the endomorphisms $d \tau_i$ on $Lie(A_K/K)$ can be
represented by the matrices $M(\tau_i)$. Since all eigenvalues of
$d \tau_i$ on $Lie(A_K/K)$ lie in $K$ the same is true about the
eigenvalues of the matrices $M(\tau_i)$. On the other hand
Equation \ref{8347t} implies that
\[\tau^*_i(\phi^e (L_j))=\sum_{s=1}^g \phi^e(m_{js}(\tau_i))
\phi^e(L_s)\] for $e=1,2$. Hence $\tau_i$ act on the space $V$ in
Equation \ref{spnn} via a matrix of the form
\[\left[
\begin{array}{ccc}M(\tau_i) & 0 & 0\\
0 & \phi M(\tau_i) & 0\\
0 & 0 & \phi^2 M(\tau_i) \end{array} \right]\] The above matrices
have all their eigenvalues in $K$. We deduce that all the
eigenvalues of $\tau_i$ on $\bX^2(A) \otimes K$ are in $K$. On the
other hand, since each $\tau_i$ annihilates a polynomial with
$\bZ-$coefficients having simple complex roots only, the same will
be true about $\tau_i$ acting on $V$ hence the minimal polynomial
of $\tau_i$ on $V \otimes_K K^a$ has simple roots only so the
minimal polynomial of $\tau_i$ acting on $\bX^2(A) \otimes K^a$
has simple roots only. So $\tau_i$ acting on $\bX^2(A) \otimes K$
is diagonalizable over $K^a$ (by Remark \ref{tutu}, assertion 1)
hence over $K$ (by Remark \ref{tutu}, assertion 2). Then, by
Remark \ref{tutu}, assertion 3, $\tau_1,...,\tau_s$ acting on
$\bX^2(A) \otimes K$ are simultaneously diagonalizable over $K$.
So there exist at least $g$ $K-$linearly independent
$\d-$characters $\psi_1,...,\psi_g \in \bX^2(A) \otimes K$ such
that $\tau_i^* \psi_j=\lambda_{ij} \cdot \psi_j$, $\lambda_{ij}
\in K$. Multiplying $\psi_j$ by a power of $p$ we may assume
$\psi_j \in \bX^2(A)$. Now, since $End(A/R) \subset
Mat_{2g}(\bZ)$, $\tau_i$ are integral over $\bZ$ hence so are the
$\lambda_{ij}$'s. Since $R$ is integrally closed, $\lambda_{ij}
\in R$. Clearly then, for $\tau \in {\mathcal T}$, we have $\tau^*
\psi_j=\e^j(\tau) \cdot \psi_j$ for some $\e^j(\tau) \in R$ and
$\e^j$ defines a ring homomorphism ${\mathcal T} \ra R$.
\end{proof}

\section{Proof of  Theorem \ref{main}}

Fix a cusp $P^0$ of $X_1(N)_{\bC}$ and consider the morphisms
\begin{equation}
\label{morphs} X_1(N)_{\bC}  \stackrel{\beta}{\ra} J_1(N)_{\bC}
\stackrel{\pi}{\ra} A_{\bC},
\end{equation}
where $\beta$ is the Abel-Jacobi map sending $P^0 \mapsto 0$ and
$A$ is as in Theorem \ref{eichshi}. Let  $M \in \bZ$ be divisible
by $N$ such that the embedding $\iota:\cO_f \ra End(A)$, all the
cusps of $X_1(N)_{\bC}$ and all the objects and morphisms in
Equation \ref{morphs} have (compatible) models over
$\bZ[1/M,\zeta_N]$.

Let $p \in \bZ$ be any prime not dividing $M$ and unramified in
$F:=\tilde{K}_f(\zeta_N)$ where $\tilde{K}_f$ is the normal
closure of $K_f$ in $\bC$. Assume we are given an embedding
$\rho:\cO_f[1/N,\zeta_N] \ra R_p=\hat{\bZ}_p^{ur}$. The latter can
be lifted to an embedding $\tilde{\rho}:\cO_F[1/M] \ra R_p$ where
$\cO_F$ is the ring of integers of $F$. Then the morphisms in
Equation \ref{morphs} induce (by base change via
$\tilde{\rho}$)
morphisms of $R_p-$schemes
\begin{equation}
\label{2} X_1(N)_{R_p} \stackrel{\beta}{\ra} J_1(N)_{R_p}
\stackrel{\pi}{\ra}  A_{R_p}.
\end{equation}
Select primes $l_1,...,l_s$ not dividing $Np$ such that the
endomorphisms \[T(l_1)_*,...,T(l_s)_*\] of $J_1(N)/\bZ[1/N]$
generate the subring
\[\bZ[T(l')_*\ | \ l'\ \  prime \ \ not\ \  dividing \ \ Np]\]
 of
$End(J_1(N)/\bZ[1/N])$. Let ${\mathcal T}^{(Np)}$ be the subring
of $End(A_{R_p}/R_p)$ generated by  $\iota(a_{l'})$ with $l'$
prime not dividing $Np$. Clearly ${\mathcal T}^{(Np)}$ is
generated as a ring by $\iota(a_{l_1}),..., \iota(a_{l_s})$.

{\bf Claim.} {\it  Each of the maps
$d(\iota(a_{l_1})),...,d(\iota(a_{l_s})) \in
End(Lie(A_{K_p}/K_p))$ has all its eigenvalues in $K_p$, the
fraction field of $R_p$.}

 Indeed an eigenvalue of
$d(\iota(a_{l_i}))$ on $Lie(A_{K_p}/K_p)$ is also an eigenvalue of
$d(\iota(a_{l_i}))$ on $Lie(A_{\bQ}/\bQ)$ hence on
$H^0(A_{\bC},\Omega)$. By the Eichler-Simura Theorem \ref{eichshi}
the latter identifies with the $\bC-$linear span $V_f$ of
$\{f^{\sigma}\ |\ \sigma:K_f \ra \bC\}$ and the action of
$d(\iota(a_{l_i}))$ corresponds to the action on $V_f$ of
$T_2(l_i)$. But the eigenvalues of $T_2(l_i)$ on $V_f$ are clearly
in $\tilde{K}_f$ and our Claim is proved.

Now each $a_{l_i}$, being an element of the field $K_f$,
annihilates a polynomial with $\bZ-$coefficients having simple
complex roots only. Hence the same is true about $\iota(a_{l_i})$.
 By Proposition
\ref{eigch} there exist $g=[K_f:\bQ]$
 $R_p-$linearly independent  $\d-$characters of order $2$,
\begin{equation}
\label{3} \psi_1,...,\psi_g:A_{R_p} \ra {\bf
G}_{a,R_p}=\bA^1_{R_p}
\end{equation}
and  ring homomorphisms
$\e^1,...,\e^g:{\mathcal T}^{(Np)} \ra R_p$
 such that
for any prime $l'$ not dividing $Np$ we have
\[\psi_j \circ
\iota(a_{l'})=\e^j(\iota(a_{l'})) \cdot \psi_j,\ \ j=1,...,g.\]
 We may then consider the
$\d-$morphisms of order $2$,
\begin{equation}
\label{4} f_j^{\sharp}:=\psi_j \circ \pi  \circ \beta:X_1(N)_{R_p}
\ra \bA^1_{R_p},\ \ \ j=1,...,g.
\end{equation}
Their restrictions to $Y_1(N)_{R_p}$ can be viewed as $\d-$modular
forms of weight $0$.

 Let now $F'$ be a number field containing $F$, let $v'$ be a
valuation on $F'$ above $p$ and let $P$ be any $\cO_{F',v'}-$point
of $Y_1(N)$ where $\cO_{F',v'}$ is the valuation ring  of $F'$ at
$v'$. Since $\cO_F[1/Ml] \subset \cO_{F',v'}$
 all cusps of $X_1(N)$ are $\cO_{F',v'}-$points and
$T(l)_*$ has a model over $\cO_{F',v'}$. So one may write
\[\sigma_{2*}\sigma_1^* P=\sum_{i=1}^{l+1} P_i,\ \ \
 \sigma_{2*}\sigma_1^* P^0=\sum_{i=1}^{l+1} P_i^0,\] where $P_i$
are distinct $\cO_{F'',v''}-$points of $X_1(N)$, $\cO_{F'',v''}$
an unramified finite extension of $\cO_{F',v'}$, and $P_i^0$ are
(a priori not necessarily distinct) cusps of $X_1(N)$ so they are
$\cO_{F',v'}-$points. So, for $l$ prime, not dividing $Np$,
Equation \ref{forforT} yields
\begin{equation}
\label{1}
T(l)_*(\beta(P))=\sum_{i=1}^{l+1}\beta(P_i)-\sum_{i=1}^{l+1}\beta(P_i^0).
\end{equation}
 Choose now an  embedding
 $\cO_{F'',v''} \ra R_p$ extending the fixed embedding
 $\cO_F[1/M] \ra R_p$ above;
in particular the points $P_i$ and $P^0_i$ can be viewed as
 $R_p-$points of $X_1(N)$.
By Theorem \ref{mandrin}, $\beta(P_i^0)$ are torsion points hence
their image via $\psi_j \circ \pi$ is $0$. Using this and the fact
that $\psi_j$ and $\pi$ are homomorphisms, we get, by Equation
\ref{1}:
\begin{equation}
\label{88}
\begin{array}{rcl}
\psi_j \pi T(l)_*(\beta(P)) & = &  \psi_j \pi
(\sum_{i=1}^{l+1}\beta(P_i)-\sum_{i=1}^{l+1}\beta(P_i^0))\\
\  & \  & \ \\
\  & = & \sum_{i=1}^{l+1} f_j^{\sharp} (P_i)\\
\  & \  & \  \\
\  & = & (T(l)f_j^{\sharp})(P).\end{array}\end{equation}
 On the other hand, using Theorem
\ref{eichshi}, we get
\begin{equation}
\label{77} \psi_j \pi T(l)_*(\beta(P))= \psi_j(\iota(a_l) \pi
\beta (P))= \e^j(\iota(a_l)) \cdot f_j^{\sharp}(P).
\end{equation}
Equations \ref{77} and \ref{88} imply that
\begin{equation}
\label{ce} (T(l)f_j^{\sharp})(P)= \e^j(\iota(a_l)) \cdot
f_j^{\sharp}(P).
\end{equation}
In the above equality
 $P$ has coordinates in $\cO_{F',v'}$. Since $F'$ and $v'$ are
 arbitrary, and since, by smoothness,
 the image of
 \[\bigcup_{F',v'} Y_1(N)(\cO_{F',v'}) \ra Y_1(N)(R_p)\]
  is $p-$adically dense in $Y_1(N)(R_p)$
 it  follows, by continuity, that Equation \ref{ce}  holds for any
 $R_p-$point $P$ of $Y_1(N)$. This implies that
 $T(l)f_j^{\sharp}=\e_j(a_l)
 \cdot f_j^{\sharp}$ where $\e_j:\cO_f^{(Np)} \ra R_p$ is the
 composition
 \[\cO_f^{(Np)} \stackrel{\iota}{\ra} {\mathcal T}^{(Np)}
 \stackrel{\e^j}{\ra} R_p.\]
 To conclude it is enough to check that
 $f^{\sharp}_1,...,f^{\sharp}_g$ are $R_p-$linearly independent.
 Assume $\sum c_jf^{\sharp}_j=0$, $c_j \in R_p$.
For large $n$, the image of the
 natural map
 \[X_1(N)(R_p)^n \ra J_1(N)(R_p)\]
 contains all the $R-$points of an open subset $U$ of $J_1(N)$.
 It follows that the restriction of $(\sum c_j \psi_j) \circ \pi$
 to $U(R_p)$ is $0$. This immediately implies that $(\sum c_j \psi_j) \circ
 \pi=0$. Hence $\sum c_j \psi_j=0$ which implies $c_j=0$ for all
 $j$. This ends our proof.

\begin{remark}
\label{ccuusspp} Note that the $\d-$modular forms $f_j^{\sharp}$
constructed above are $\d-$cusp forms (with respect to $\rho$) in
the sense of Definition \ref{cussp}. Indeed $f^{\sharp}_j$ come
from $\d-$ morphisms
\[\psi \circ \pi \circ \beta:X_1(N)_{R_p} \ra
\bA^1_{R_p},\] so they are $\d-$holomorphic. Now,
  by the
Manin-Drinfeld Theorem \ref{mandrin}, the images of the cusps in
$J_1(N)(R)$ via the Abel-Jacobi map
\[\beta:X_1(N)(R_p) \ra J_1(N)(R_p)\]
are torsion points. On the other hand
\[\psi \circ \pi:J_1(N)(R_p) \ra {\bf G}_{a,R_p}(R_p)=R_p\]
 is a homomorphism with torsion free target so it vanishes on
 torsion points.
\end{remark}

\section{$\d-$Fourier expansions}

We start by recalling the background of {\it classical} Fourier
expansions; cf. \cite{DI}, p. 112. (The discussion in loc. cit.
involves the model $X_{\mu}(N)$ instead of the model $X_1(N)$ used
here but the two models, and hence the two theories, are
isomorphic over $\ZN$ cf. \cite{DI}, p. 113.) There is a point
$s_{\infty}:\ZN \ra X_1(N)_{\ZN}$ arising from the generalized
elliptic curve $\bP^1_{\ZN}$ with its canonical embedding of
$\mu_{N,\ZN} \simeq (\bZ/N\bZ)_{\ZN}$; the complex point
corresponding to $s_{\infty}$ is the cusp $\Gamma_1(N) \cdot
\infty$. The map $s_{\infty}$ is a closed immersion; denote by
$\tilde{X}_1(N)_{\ZN}$ the completion of $X_1(N)_{\ZN}$ along the
image of $s_{\infty}$. Now consider the Tate generalized elliptic
curve
\[Tate(q)/\ZN [[q]];\]
it has a canonical immersion $\alpha_{can}$ of $\mu_{N,\ZN} \simeq
(\bZ/N\bZ)_{\ZN}$  so there is an induced map $Spec\ \ZN \ra
X_1(N)_{\ZN}$ (which further composed with $Spec\ \ZN \ra Spec\
\ZN [[q]]$, $q \mapsto 0$, equals $s_{\infty}$). There is an
induced isomorphism
\begin{equation}
\label{ttt} Spf\ \ZN[[q]] \ra \tilde{X}_1(N)_{\ZN}.
\end{equation}
There is a canonical $1-$form $\omega_{can}$ on the elliptic curve
$Tate(q)/\ZN ((q))$ over $\ZN ((q)):=\ZN [[q]][1/q]$ such that the
induced map
\[S_2(\Gamma_1(N),\bC) \ra \bC((q)), \ f \mapsto
f_{\infty}:=f_{\infty}(q):= f(Tate(q)/\bC((q)),\alpha_{can},
\omega_{can}),\] has image in $q\bC[[q]]$ and is the classical
Fourier expansion at the cusp $\Gamma_1(N) \cdot \infty$. (Here we
interpret $f$ as a function of triples in the style of
\cite{Katz}.)

Next we move to the $\d-$theory. Fix a prime $p$ not dividing $N$
and fix a homomorphism $\rho:\ZN \ra R_p$. We may consider the
prolongation sequence $S_{\infty}^*$ defined by
\[S_{\infty}^r:=R_p((q))\hat{\ }[q',q'',...,q^{(r)}]\hat{\ },\]
where $\hat{\ }$ means $p-$adic completion and
$q',q'',...,q^{(r)}$ are new indeterminates; $S_{\infty}^{r+1}$ is
viewed as an $S_{\infty}^r-$algebra via the inclusion and the
$p-$derivations $\d:S_{\infty}^r \ra S_{\infty}^{r+1}$ are the
unique $p-$derivations satisfying $\d q^{(i)}=q^{(i+1)}$.
Explicitly, they are defined as follows. First extend $\phi:R_p
\ra R_p$ to ring homomorphisms $\phi:S_{\infty}^r \ra
S_{\infty}^{r+1}$ by requiring that
\begin{equation}
\label{nitzel} \phi(q):=q^p+pq',\ \phi(q'):=(q')^p+pq'',...
\end{equation}
 and then define $\d:S_{\infty}^r \ra S_{\infty}^{r+1}$ by
 \begin{equation}
 \label{incanitzel}
\d F:= \frac{\phi(F)-F^p}{p}. \end{equation}
 Finally  define the
{\it $\d-$Fourier expansion map}
\[E_{\infty,\rho}:M^r(\Gamma_1(N),R_p,w) \ra S_{\infty}^r,\ \ \ f \mapsto
E_{\infty,\rho}(f)=f_{\infty,\rho},\] by the formula
\[f_{\infty,\rho}:=f_{\infty,\rho}(q,q',...,q^{(r)}):=
f(Tate(q)/S_{for}^0,\alpha_{can},\omega_{can},S_{\infty}^*).\]
Note that if $f \in S_m(\Gamma_1(N),\bC)$ is a newform and
$\rho:\cO_f[1/N,\zeta_N] \ra R_p$ extends the
homomorphism  $\rho:\ZN \ra R_p$
above  then, by functoriality,
$(f^{\rho})_{\infty,\rho}(q)=(f_{\infty}(q))^{\rho}$,
where $(f^{\rho})_{\infty,\rho}(q)$ is the $\d-$Fourier expansion
of $f^{\rho} \in M^0(\Gamma_1(N),R_p,m)$ and
$(f_{\infty}(q))^{\rho}$ is the image of the classical Fourier
expansion $f_{\infty}(q) \in \cO_f[[q]]$ via the natural map
$\rho:\cO_f[[q]] \ra R_p[[q]]$. If $\cO_f=\bZ$ we simply have
$(f_{\infty}(q))^{\rho}=f_{\infty}(q)$.
Also note that the $\d-$Fourier expansion map commutes with $\phi$
i.e.
$(f^{\phi})_{\infty,\rho}=(f_{\infty,\rho})^{\phi}$.

\begin{lemma}
\label{deltaexp} ($\d-$expansion principle) The $\d-$Fourier
expansion map
\[E_{\infty,\rho}:M^r(\Gamma_1(N),R_p,w) \ra S^r_{\infty}\]
is
injective, with torsion free cokernel.
\end{lemma}

\begin{proof}
This
follows from \cite{book}, Proposition 4.43, by ``taking fractions"
with denominators powers of a local equation defining the cusp
$\infty$.
\end{proof}

For two elements $F,G \in S^r_{\infty}$ we write $F \sim G$ if
$F=\lambda \cdot G$ with $\lambda \in R_p^{\times}$ and we write
$\bar{F} \sim \bar{G}$ if the images,
\[\bar{F}, \bar{G} \in S^r_{\infty} \otimes k=
k((q))[q',q'',...,q^{(r)}]\] of $F$ and $G$ satisfy $\bar{F}=c
\cdot \bar{G}$ for some $c \in k^{\times}$, $k:=R_p/pR_p$. Let
\[\Psi=\Psi(q,q'):=\frac{1}{p}
\cdot log\left( 1+p\frac{q'}{q^p} \right)=\sum_{n=1}^{\infty}
(-1)^{n-1} p^{n-1} n^{-1} \left( \frac{q'}{q^p} \right)^n \in
R_p((q))\hat{\ }[q']\hat{\ }.\]
 Then the $\d-$Fourier expansion of $f^r =
f^r_{crys} \in M^r(\Gamma_1(N),R_p,-1-\phi^r)$ is given by:

\begin{lemma}
\cite{Barcau}
 \label{olema} \[f^r_{\infty,\rho} \sim
 \Psi^{\phi^{r-1}}+p\Psi^{\phi^{r-2}}+...+p^{r-1} \Psi  \in
R_p((q))\hat{\ }[q',...,q^{(r)}]\hat{\ }.\] In particular
\[\overline{f^r_{\infty,\rho}} \sim \left( \frac{q'}{q^p}\right)^{p^{r-1}} \in
k((q))[q',...,q^{(r)}].\]
\end{lemma}

Next we would like to compute the $\d-$Fourier
expansion of the forms $f^{\sharp}=f^{\sharp}_1$
in Theorem \ref{main} for $g=1$.

\begin{theorem}
\label{fur}
 Let $f =\sum a_n q^n \in S_2(\Gamma_0(N),\bC)$ be a newform
with $K_f=\bQ$ which is not of CM type. Then, for any sufficiently
large prime $p$ and any embedding $\rho:\bZ[1/N,\zeta_N] \ra R_p$
there is a unique $\d-$eigenform $f^{\sharp} \in
M^2(\Gamma_1(N),R_p,0)$  with $\d-$Fourier expansion:
\begin{equation} \label{h1}
f^{\sharp}_{\infty,\rho}(q,q',q'') = \frac{1}{p} \sum_{n\geq 1}
\frac{a_n}{n} (q^{n\phi^2} -a_p q^{n\phi} +p q^n) \in
 R_p((q))\hat{\ }[q',q'']\hat{\ }.
 \end{equation}
 Moreover $f^{\sharp}$ is a $\d-$cusp form (with respect to
 $\rho$)
 belonging to $f$ with exponent $0$ and
 \begin{equation}
 \label{h12}
f^{\sharp}_{\infty,\rho}(q,0,0)=\sum_{(n,p)=1} \frac{a_n}{n} q^n
\in R_p[[q]].
 \end{equation}
\end{theorem}

 \begin{remark}
 1) Explicitly, in Equation \ref{h1} we have:
\[\begin{array}{rcl}
q^{n\phi} & = & (q^p+pq')^n,\\
q^{n\phi^2} & = & [(q^p+pq')^p +p(q')^p+p^2q'']^n.
\end{array}\]

2) Uniqueness of $f^{\sharp}$ follows, of course, from the
$\d-$expansion principle (Lemma \ref{deltaexp}).

3) The series
\begin{equation}
\label{nocnoc} \sum_{(n,p)=1} \frac{a_n}{n} q^n
\end{equation}
 is normalized and has
coefficients in $\bZ_{(p)}$ but not all its coefficients are in
$\bZ$. Indeed if the latter were the case then any prime $l \neq
p$ would divide $a_l$. Since, for big enough $l$, $a_l$ are the
traces of Frobenius of an elliptic curve $A$ over $\bQ$ taken
modulo $l$ it would follow that $A$ has supersingular reduction
for sufficiently big $l$, a contradiction. In particular the
series \ref{nocnoc} is not a (classical) eigenform.
\end{remark}

\begin{proof}
We begin with a preparatory discussion; in this discussion we will
not assume yet that $f$ is not of CM type (for we will use this
discussion later in case $f$ is of CM type). We place ourselves in
the context of the proof of Theorem \ref{main}. Since $K_f=\bQ$
the field $F$ in that proof equals $\bQ(\zeta_N)$. Also we
 let $\beta:X_1(N)_{\bC} \ra J_1(N)_{\bC}$ be the Abel-Jacobi map that sends
the cusp $P^0=\infty$ into $0$. One can choose $A$ in Theorem
\ref{eichshi}, and hence in the proof of Theorem \ref{main},
such
that
\[\Phi:=\pi \circ \beta:X_1(N)_{\bC} \ra A_{\bC}\]
satisfies
\[\Phi^* \omega_A=c \cdot 2 \pi i \cdot f(z)dz=
c \cdot \sum_{n \geq 1} a_n q^{n-1} dq,\] where $\omega_A$ is a
$1-$form on $A$ over $\bQ$, $q=e^{2 \pi i z}$ and $c \in
\bQ^{\times}$. Cf. \cite{darmon}, p. 19. Let $T$ be an \'{e}tale
coordinate around the origin $0$ of $A$ such that $T$ vanishes at
$0$. Let $L(T) \in \bQ[[T]]$ be the logarithm of the formal group
of $A$ associated to $T$; cf. \cite{sil}. Then
$\omega_A=\omega(T)dT \in \bQ[[T]] \cdot dT$; replacing $\omega_A$
by a $\bQ^{\times}-$multiple of it we may (and will) assume
$\omega(0)=c$. Let $m \in \bZ$ be such that $X_1(N)$, $A$, $T$,
$\omega_A$ have (compatible)
 models over $\bZ[1/m]$. Then
 $P^0$ and $\Phi$ are defined over $\cO:=\bZ [1/Nm,\zeta_N]$.
  We have an induced homomorphism
$\Phi^*:\cO[[T]] \ra \cO[[q]]$
and we set
$\Phi^*(T)=\varphi(q) \in \cO[[q]]$.
Since $\omega_A=c \cdot \frac{dL}{dT} \cdot dT$ we get
\[c \cdot (\sum_{n \geq 1} a_n q^{n-1}) dq=
\Phi^* \omega_A=c \cdot \frac{dL}{dT}(\varphi(q)) \cdot \frac{d
\varphi}{dq} (q) \cdot dq= c \cdot \frac{d}{dq}(L(\varphi(q))).\]
Setting $q=0$ in the coefficients of $dq$  we get that
$\frac{d \varphi}{dq}(0)=1$.
Also we deduce that
\[L(\varphi(q))=
\sum_{n \geq 1} \frac{a_n}{n} q^n.\] Now, for $p$ sufficiently
large,  $\rho:\ZN \ra R_p$  induces  a homomorphism $\rho:
\bZ[1/Nm,\zeta_N] \ra R_p$.

From this point on we assume $f$ is not of CM type.
Since $A$ is
not a CM elliptic curve it follows that $A_{R_p}$ does not have a
lift of Frobenius (i.e. there is no morphism of schemes
$\tilde{\phi}: A \ra A$ over $\bZ$ lifting the morphism $Spec\ R_p
\ra Spec\ R_p$ induced by $\phi$ such that the reduction mod $p$
of $\tilde{\phi}$ is the $p-$power Frobenius on $A_{R_p} \otimes
k$.) Since $A_{R_p}$ does not have a lift of Frobenius,  by
\cite{book}, Theorem 7.22 and \cite{cr97}, Theorem 1.10, one can
assume the $\d-$character $\psi$ in the proof of Theorem
\ref{main} gives rise to the series (still denoted by)
\[\psi=\frac{1}{p} (\phi^2-a_p\phi+p) L(T) \in R_p[[T]][T',T'']\hat{\ },\]
where $\phi$ is viewed here as naturally extended to series with
$K-$coefficients. Then the $\d-$Fourier expansion of the form
$f^{\sharp}:=f^{\sharp}_1$ provided by the proof of Theorem
\ref{main} equals
\begin{equation}
\label{titi}
\begin{array}{rcl}
 f^{\sharp}_{\infty,\rho} & = & \Phi^{*}
\psi=\Phi^{*}\left(\frac{1}{p} (\phi^2-a_p\phi+p) L(T)\right)
\\
\  &  \ & \  \\
\  & = &
\frac{1}{p}
\{(\phi^2-a_p\phi+p)(L(T))\}(\varphi(q),\d(\varphi(q)),\d^2(\varphi(q)))
\\
\  &  \ & \  \\
\  & = & \frac{1}{p} (\phi^2-a_p\phi+p) (L(\varphi(q)))
= \frac{1}{p}  (\phi^2-a_p\phi+p) \left( \sum_{n\geq 1}
\frac{a_n}{n} q^n \right)\\
\  &  \ & \  \\
 \  & = & \frac{1}{p} \sum_{n\geq 1}
\frac{a_n}{n} (q^{n\phi^2} -a_p q^{n\phi} +p q^n)\end{array}
\end{equation} and our
Equation \ref{h1} follows.

Setting $q'=q''=0$ in this equation we get
\[f^{\sharp}_{\infty,\rho}(q,0,0)=\frac{1}{p} \sum_{n \geq 1}
\frac{a_n}{n} (q^{np^2}-a_pq^{np}+pq^n)=\sum_{c \geq 1} c_nq^n,\]
\[c_n=\frac{1}{p} \left( \frac{a_{n/p^2}}{n/p^2}-
\frac{a_{n/p}a_p}{n/p}+\frac{a_n}{n} p \right).\] Here, by
definition,
 $a_x=0$ for $x \in \bQ \backslash \bZ$. Note that for $p$ not
 dividing $n$ we get $c_n=a_n/n$. For $n=pm$ with $p$ not dividing
 $m$ we get
 \[c_n=\frac{a_n}{n}-\frac{a_ma_p}{pm}=0.\]
 For $n=p^im$ with $p$ not dividing $m$ and $i \geq 2$ we get,
by Equation \ref{tzu}, that
 \[c_n=\frac{a_{p^{i-2}}a_m}{p^{i-1}m}-
 \frac{a_{p^{i-1}m} a_p}{n} +\frac{a_n}{n}=
 \frac{a_m}{n}(pa_{p^{i-2}}-a_{p^{i-1}}a_p+a_{p^i})=0.\]
 Equation \ref{h12} follows.
\end{proof}

In the CM case we have the following:

\begin{theorem}
\label{fur2}
 Let $f =\sum a_n q^n \in S_2(\Gamma_0(N),\bC)$ be a newform
with $K_f=\bQ$ which is  of CM type (corresponding to an imaginary
quadratic field ${\mathcal K}$). Then, for any sufficiently large
prime $p$ and any embedding $\rho:\bZ[1/N,\zeta_N] \ra R_p$ there
is a unique $\d-$eigenform $f^{\sharp} \in M^2(\Gamma_1(N),R_p,0)$
such that the following hold:

1) If $p$ does not split  in ${\mathcal K}$ then Equation \ref{h1}
holds (with $a_p=0$).

2) If $p$ splits in ${\mathcal K}$ and if $pu$ is the unique root
in $p\bZ^{\times}_p$ of the polynomial $x^2-a_px+p$ then the
following equation holds:
\begin{equation} \label{H1}
f^{\sharp}_{\infty,\rho}(q,q',q'') = \frac{1}{p} \sum_{n\geq 1}
\frac{a_n}{n} (q^{n\phi} -pu q^n) \in
 R_p((q))\hat{\ }[q',q'']\hat{\ }.
 \end{equation}
 Moreover $f^{\sharp}$ is a $\d-$cusp form (with respect to
  $\rho$), $f^{\sharp}$ belongs to $f$ with exponent
 $0$, and the following hold:

 1') If $p$ does not split  in ${\mathcal K}$ then Equation \ref{h12}
holds;

2') If $p$ splits in ${\mathcal K}$ then
 \begin{equation}
 \label{H12}
f^{\sharp}_{\infty,\rho}(q,0,0)=-u \cdot \sum_{(m,p)=1}
\frac{a_m}{m}  \sum_{i \geq 0} u^iq^{mp^i} \in R_p[[q]].
 \end{equation}
\end{theorem}

\begin{proof}
Let us go back to the preparatory discussion in the proof of
Theorem \ref{fur}. In our case here $A$ is an elliptic curve over
$\bQ$ with CM by an order of ${\mathcal K}$. By \cite{sil2}, p.
180, there exists an elliptic curve $A'$ over $\bQ$ and an isogeny
defined over $\bQ$, $A \ra A'$, such that $A'$ has CM by the ring
of integers $\cO_{\mathcal K}$ of ${\mathcal K}$. Replacing $A$ by
 $A'$ we may assume that $A$ itself has CM by
$\cO_{\mathcal K}$. Since $A$ is defined over $\bQ$, ${\mathcal
K}$ has class number one; cf. \cite{sil2}, pp. 118, 121. By
standard facts about elliptic curves with CM (\cite{sil2}, pp 184,
133) we have that for $p$ large enough the following hold:

I) If $p$ does not split in ${\mathcal K}$ then $A_{R_p}$ has
supersingular reduction, $a_p=0$,  and $A_{R_p}$ doesn't have a
lift of Frobenius.

II) If $p$ splits in ${\mathcal K}$ then $A_{R_p}$ has ordinary
reduction, $a_p \not\equiv 0$ mod $p$, and $A_{R_p}$ has a lift of
Frobenius.

If we are in case I the argument in the proof of Theorem \ref{fur}
applies. Assume we are in case II. Then, by \cite{cr97}, Theorem
1.10, $\psi$  gives rise to the series
 (still denoted by)
\[\psi=\frac{1}{p} (\phi-up) L(T) \in R_p[[T]][T']\hat{\ }.\]
Then the $\d-$Fourier expansion of the form $f^{\sharp}$ provided
by the proof of Theorem \ref{main} equals
\begin{equation}
\label{tititi}
\begin{array}{rcl}
f^{\sharp}_{\infty,\rho} & = &  \Phi^*
\psi\\
\  & \  & \  \\
\  & = & \Phi^*\left(\frac{1}{p} (\phi-up) L(T)\right)
\\
\  & \  & \  \\
\  & = & \frac{1}{p} \{(\phi-up)(L(T))\}(\varphi(q),\d(\varphi(q)))\\
\  & \  & \  \\
\  & = & \frac{1}{p} (\phi-up) L(\varphi(q))\\
\  & \  &  \ \\
\  & = &
\frac{1}{p}  (\phi-up) \left( \sum_{n\geq 1} \frac{a_n}{n} q^n
\right)\\
\  & \  & \  \\
\  & = & \frac{1}{p} \sum_{n\geq 1} \frac{a_n}{n} (q^{n\phi} -up
q^n)\end{array}\end{equation}
 and our Equation \ref{H1} follows.

Setting $q'=q''=0$ in this equation we get
\[f^{\sharp}_{\infty,\rho}(q,0,0)=\frac{1}{p} \sum_{n \geq 1}
\frac{a_n}{n} (q^{np}-upq^n)=\sum_{c \geq 1} c_nq^n,\]
\[c_n=\frac{a_{n/p}-ua_n}{n} .\]
 Note that for $p$ not
 dividing $n$ we get
 $c_n=-ua_n/n$.
 For $n=p^im$ with $p$ not dividing $m$ and $i \geq 1$ we get,
by Equation \ref{tzu}, that
 \[c_n=\frac{a_m}{m} \cdot \frac{a_{p^{i-1}}-ua_{p^i}}{p^i}.\]
 We claim that
 \[a_{p^{i-1}}-ua_{p^i}=-p^iu^{i+1},\]
 and this will, of course, end the proof of the equality
in Equation \ref{H12}.
 To check the claim note that, for $i=1$, we get $a_1-ua_p=-pu^2$.
In general we proceed
 by induction, using Equations \ref{tzu}:
 \[\begin{array}{rcl}
a_{p^i}-ua_{p^{i+1}} & = & a_{p^i}-u(a_{p^i}a_p-pa_{p^{i-1}})\\
\  & \  & \  \\
\  & = &
 (1-ua_p)a_{p^i}+upa_{p^{i-1}}\\
 \  & \  & \  \\
 \ & = & -pu^2a_{p^i}+upa_{p^{i-1}}\\
 \  & \  & \  \\
\  & = & up(a_{p^{i-1}}-ua_{p^i})\\
\  & \  & \  \\
\  & = & -p^{i+1}u^{i+2},\end{array}\]
 and our claim is proved.
\end{proof}

\section{Independence of  $f^{\sharp}$ from $f$ and $f^{(-1)}$}

Using $\d-$Fourier expansions it is possible to prove a variant of
the   result on the independence of $f^{\sharp}$ from $f$
contained in assertion 2 of Theorem \ref{main} and also a result
on the independence of $f^{\sharp}$ from $f^{(-1)}$.

\begin{theorem}
\label{doii} Let $f \in S_2(\Gamma_0(N),\bC)$ be a newform of
weight $2$ on $\Gamma_0(N)$, $N>4$, with $K_f=\bQ$. For any
sufficiently large prime $p$ and any embedding
$\rho:\bZ[1/N,\zeta_N] \ra R_p$ let
$f^{\sharp} \in M^2(\Gamma_1(N),R_p,0)$ be the unique form
in Theorems \ref{fur} or \ref{fur2} according as $f$ is of non-CM
or CM type respectively.
Then the following hold

1)  There is no $G \in \cJ$ such that $G \cdot
 f^{\sharp}$ belongs to the $J-$linear span of
\[\{f^{\rho},f^{\rho \phi},f^{\rho\phi^2},f^{\rho \phi^3},...\}.\]

2) $f^{\sharp}$ does not belong to the $I_{ord}-$linear span in
$M^r_{ord}(\Gamma_1(N),R_p,0)$ of
\[\{f^{(-1)\rho},f^{(-1)\rho \phi},
f^{(-1)\rho\phi^2},f^{(-1)\rho \phi^3},...\}.\]
\end{theorem}

\begin{proof}
We prove assertion 1.
 Since there are infinitely many
primes $l$ of ordinary reduction for the elliptic curve $A$ we may
choose two such distinct primes $l_1$ and $l_2$; so
$a_{l_1}a_{l_2} \neq 0$. Let $p$ be a prime which is sufficiently
big so  that Theorem \ref{fur} holds in case $A$ doesn't have  CM
(respectively Theorem \ref{fur2} holds in case $A$ has CM) and, in
addition, $p$ does not divide $l_1l_2(l_1-l_2)a_{l_1}a_{l_2}$. Fix
a homomorphism
 $\rho:\bZ[1/N,\zeta_N] \ra R_p$ and let $f^{\sharp}$ be the unique form
 in  Theorem \ref{fur} in case $A$ doesn't have CM (respectively
 the unique form in
Theorem \ref{fur2} in case $A$ has CM).

Assume $A$ doesn't have CM; the case when $A$ has CM is entirely
similar and left to reader. Assume
there exists $G \in \cJ$ such that $G \cdot f^{\sharp}$ is a
$J-$linear combination of forms $f^{\rho \phi^d}$. Hence
$G \cdot f^{\sharp}$ is
 a
$R_p-$linear combination of forms $F \cdot f^{\rho\phi^d}$ with $F
\in \cJ$. By taking $\d-$Fourier expansions we get that
$G_{\infty,\rho} \cdot f^{\sharp}_{\infty,\rho}$ is an
$R_p-$linear combination of forms $F_{\infty,\rho} \cdot
f^{\rho\phi^d}_{\infty,\rho}$. Reducing modulo $p$, setting
$q''=0$, and using
 Lemma
\ref{olema} and Theorem \ref{fur} we have a
congruence mod $p$ of the form
\[ \left( \frac{q'}{q^p} \right)^s \left( \sum_{(m,p)=1}
\frac{a_m}{m} q^m+q'B(q,q') \right) \equiv  \sum_{j,d \geq 0}
\lambda_{dj} \left( \frac{q'}{q^p} \right)^j \left( \sum_{m \geq
1} a_m q^{mp^d} \right)
\]
in $R_p((q))[q']\hat{\ }$, where $s \geq 0$, $\lambda_{ij} \in
R_p$, and $B(q,q') \in R_p[[q]][q']$. Let $l$ be either $l_1$ or
$l_2$. Identifying the coefficients of $(q')^sq^{l-ps}$ in the
above Equation we get
\[\frac{a_l}{l} \equiv  \lambda_{0s} \cdot a_l\ \ \ mod\ \ (p)\]
in $R_p$. Since $a_{l} \not\equiv 0$ mod $p$ we get
\[1 \equiv l \cdot  \lambda_{0s}\ \ \
mod\ \ (p).\]
So $l_1 \equiv l_2$ mod $p$, a contradiction.

We prove assertion 2. Assume $f^{\sharp}=\sum_{j \geq 0}
\lambda_j f^{(-1)\rho\phi^j}$, $\lambda_j \in I_{ord}$.
By Lemma \ref{mabatelakap}
we may assume $\lambda_j \in R$. Looking
at $\d-$Fourier expansions we get one of the following equalities:
\begin{equation}
\label{proscut}
\begin{array}{rcl}
\frac{1}{p} \sum_{n \geq 1}
\frac{a_n}{n}(q^{n\phi^2}-a_pq^{n\phi}+pq^n)& = & \sum_{j \geq 0}
\sum_{(n,p)=1} \lambda_j \frac{a_n}{n} q^{n\phi^j}\\
\  & \  & \ \\
\frac{1}{p} \sum_{n \geq 1}
\frac{a_n}{n}(q^{n\phi}-puq^n) & = & \sum_{j \geq 0}
\sum_{(n,p)=1} \lambda_j \frac{a_n}{n} q^{n\phi^j}.
\end{array}
\end{equation}
Setting $q'=q''=...=0$ and picking out the coefficient of $q^{p^j}$
we get $\lambda_0=1$, $\lambda_1=\lambda_2=...=0$ in the first case and
$\lambda_j=-u^{j+1}$ in the second case respectively. In both situations
we clearly get a contradiction.
\end{proof}

\section{$\d-$Serre operators}

Recall from \cite{book}, p. 255, that the Serre-Katz operators on
modular forms \cite{Katz}, p. 169, induce $R-$derivations
\[\partial_j:M^r(\Gamma_1(N),R,*) \ra M^r(\Gamma_1(N),R,*),\ \
j \geq 0,\] such that if $w=\sum a_i\phi^i$ then
\[\partial_j M^r(\Gamma_1(N),R,w) \subset
M^r(\Gamma_1(N),R,w+2\phi^j).\]
 According to \cite{Katz} (or \cite{book}, p.
255) the {\it Ramanujan form} defines an element
\[P \in M^0_{ord}(\Gamma_1(N),R,2).\]
(N.B: the $P$ in \cite{Katz} is $12$ times the $P$ in \cite{book};
here we are using the $P$ in \cite{book}.) One can consider the
$R-$derivations
\[\partial_j^*:M^r_{ord}(\Gamma_1(N),R,*):=\bigoplus_{w \in
W(r)}M^r_{ord}(\Gamma_1(N),R,w) \ra M^r_{ord}(\Gamma_1(N),R,*),\]
where, for $w=\sum a_i \phi^i$, the restriction of $\partial_j^*$
to $M^r_{ord}(\Gamma_1(N),R,w)$ equals
\[\partial_j+a_jp^j
P^{\phi^j}.\] Recall from \cite{book}, p. 93, that one also
defines
$\partial_{**}:=\sum_{j \geq 0} p^{-j}\partial_j$.
On the other hand one can consider the $R-$derivation $\theta:=q
\frac{d}{dq}$ on $S^0_{\infty}:=R((q))\h$. Exactly as in
\cite{book}, p. 113, there exist unique $R-$derivations
\[\theta_j:\bigcup_{r\geq 0} S^r_{\infty} \ra \bigcup_{r\geq 0}
S^r_{\infty}\]
 satisfying the properties
\begin{equation}
\label{tzitzu} \theta_j \circ \phi^s=0,\ \ on \ \ S^0_{\infty}\ \
for\ \  s \neq j,
\end{equation}
\[\theta_j \circ \phi^j=p^j \cdot \phi^j \circ \theta\ \ \
\ \ on \ \ S^0_{\infty}.\]

\begin{lemma}
\label{muu} For all $j \geq 0$, $r \geq 0$, $w \in W(r)$, we have
an equality of maps
\[E_{\infty,\rho} \circ \partial_j^*=\theta_j \circ
E_{\infty,\rho}:M^r(\Gamma_1(N),R,w) \ra S^r_{\infty}.\] In
particular we have an equality of maps
\[E_{\infty,\rho} \circ \partial_j=\theta_j \circ
E_{\infty,\rho}:M^r(\Gamma_1(N),R,0) \ra S^r_{\infty}.\]
\end{lemma}

\begin{proof}
Same argument as in \cite{book}, p. 259, where the case of
Serre-Tate expansions (rather than Fourier expansions) was
considered; to make that argument work one uses \cite{Katz}, p.
180.
\end{proof}

\begin{remark}
\label{miinnc} By \cite{book}, p. 113, $\theta_j$ sends each
$R[[q]][q',...,q^{(r)}]\h$ into itself hence induces a
$K-$derivation (still denoted by) $\theta_j$ on
$K[[q,...,q^{(r)}]]$ which still satisfies Equations \ref{tzitzu}
(with $S_{\infty}^0$ replaced by $K[[q]]$).
\end{remark}

In the next two Propositions
we compute the effect of $\partial_j$ on
$f^{\sharp}$.

\begin{proposition}
\label{se1}
 Assume the hypotheses and notation of Theorem
\ref{fur}. Then
\begin{equation}
\label{soan}
\partial_2f^{\sharp}=pf^{\phi^2},\ \ \
\partial_1f^{\sharp}=-a_pf^{\phi},\ \ \ \partial_0f^{\sharp}=f.
\end{equation}
In particular
$\partial_{**}f^{\sharp}=\frac{1}{p}(f^{\phi^2}-a_pf^{\phi}+pf)$.
\end{proposition}

\begin{proof}
By Lemma \ref{muu}, Theorem \ref{fur}, and Remark \ref{miinnc} one
has:
\[\begin{array}{rcl}
(\partial_2 f^{\sharp})_{\infty,\rho} & = &
\theta_2(f^{\sharp}_{\infty,\rho})\\
\  & \  & \  \\
\  & = &
\frac{1}{p} \sum_{n \geq 1} \frac{a_n}{n} \theta_2(\phi^2(q^n))\\
\ & \  & \  \\
\  & = & \frac{1}{p} \sum_{n \geq 1} \frac{a_n}{n} p^2
\phi^2(\theta(q^n))\\
\  & \  & \  \\
\  & = & p \left( \sum_{n \geq 1} a_nq^n
\right)^{\phi^2}\\
\  & \  & \  \\
\  & = & (pf^{\phi^2})_{\infty,\rho}.\end{array}\] By the
$\d-$expansion principle (Lemma \ref{deltaexp}) we get
$\partial_2f^{\sharp}=pf^{\phi^2}$. The other equalities are
obtained in the same way.
\end{proof}

In a similar way one proves:

\begin{proposition}
\label{se2} Assume the hypotheses and notation of Theorem
\ref{fur2}. Then the following hold.

1) If $p$ splits in ${\mathcal K}$ then Equations \ref{soan} hold
(with $a_p=0$).

2) If $p$ does not split in ${\mathcal K}$ then
\begin{equation}
\label{opajj}
\partial_1f^{\sharp}=f^{\phi},\ \ \ \partial_0f^{\sharp}=-uf.
\end{equation}
In particular
$\partial_{**}f^{\sharp}=\frac{1}{p}(f^{\phi}-puf)$.
\end{proposition}

Note that
Equations \ref{soan} and \ref{opajj} together with the condition
that $f^{\sharp}$ is in $M^r(\Gamma_1(N),R_p,0)$ ($r=1,2$) pin down
$f^{\sharp}$ up to an additive constant in $R$.

\section{The Hecke operator $T(p)_{\infty}$}

A direct attempt to define the Hecke operator $T(p)$ on
$\d-$modular forms along the lines of Equation \ref{tofl}
obviously fails. The ``expected"  definition of $T(p)$ on
arbitrary series in $R[[q,q',...,q^{(r)}]]$ is also easily seen to
fail. We will define $T(p)_{\infty}$ on a certain $R-$submodule of
$R[[q,q',...,q^{(r)}]]$; then $f_{\infty,\rho}^{\sharp}$ will be
in that submodule and will turn out to be an eigenvector for
$T(p)_{\infty}$ with eigenvalue $a_p(f)$. By considering
 a slightly different
$R-$module of series we will show that the $\d-$Fourier expansions
$f^r_{\infty,\rho}$ of $f^r=f^r_{crys}$ are also eigenvectors of
an appropriate version, $T(p)_{\infty,2}$, of $T(p)_{\infty}$ with
eigenvalues $p(p+1)$.

\begin{definition}
Let $q_1,...,q_m$ be variables and let $S_1,...,S_m$ be the
fundamental  symmetric polynomials in $q_1,...,q_m$; so
\[S_1=q_1+...+q_m,...,S_m=q_1...q_m.\]
A series
\[G \in R[[q_1,...,q_m,...,q_1^{(r)},...,q_m^{(r)}]]\]
is {\it $\d-$symmetric} if there exists a series
\[G_{(m)} \in R[[q_1,...,q_m,...,q_1^{(r)},...,q_m^{(r)}]]\]
such that
\begin{equation}
\label{rupe}
G(q_1,...,q_m,...,q_1^{(r)},...,q_m^{(r)})=G_{(m)}(S_1,...,S_m,...,\d^r
S_1,...,\d^r S_m).
\end{equation}
(The series $G_{(m)}$ is trivially seen to be unique; cf.
\cite{dcc}.)
\end{definition}

On the other hand, for any series $F \in R[[q,...,q^{(r)}]]$ one
can define the series
\[\Sigma_mF:=\sum_{j=1}^m F(q_i,...,q_i^{(r)})
\in R[[q_1,...,q_m,...,q_1^{(r)},...,q_m^{(r)}]].\] The series
$\Sigma_mF$ is not $\d-$symmetric in general. But there are
important examples when $\Sigma_m F$ is $\d-$symmetric; for
instance we have:

\begin{lemma}
\label{seninar} Let $\bT=(T^1,...,T^g)$ be a $g-$tuple of
variables, let $\cF \in R[[\bT_1,\bT_2]]^g$ be a formal group law,
and let $\psi \in R[[\bT,...,\bT^{(r)}]]$ be such that
\[\psi(\cF(\bT_1,\bT_2),...,\d^r\cF(\bT_1,\bT_2))=
\psi(\bT_1,...,\bT_1^{(r)})+\psi(\bT_2,...,\bT_2^{(r)})\] in the
ring
\[R[[\bT_1,\bT_2,...,\bT_1^{(r)},\bT_2^{(r)}]].\]
Let $\varphi(q) \in R[[q]]^g$ be a $g-$tuple of series and let
\[F:=\psi(\varphi(q),...,\d^r(\varphi(q))) \in
R[[q,...,q^{(r)}]].\] Then $\Sigma_mF$ is $\d-$symmetric for all
$m \geq 2$.
\end{lemma}

\begin{proof}
An easy exercise. Cf. also \cite{dcc}.
\end{proof}

\begin{corollary}
\label{ddus} If $f^{\sharp}_{\infty,\rho}$ is an in Theorems
\ref{fur} and \ref{fur2} then $\Sigma_mf^{\sharp}_{\infty,\rho}$
is $\d-$symmetric for all $m \geq 2$.
\end{corollary}

\begin{proof}
By Equations \ref{titi} and \ref{tititi}
$f^{\sharp}_{\infty,\rho}$ can be written as
\[\psi(\varphi(q),...,\d^r(\varphi(q)))\]
with $\psi$ and $\varphi$ as in Lemma \ref{seninar} and we may
conclude by Lemma \ref{seninar}.
\end{proof}

\begin{definition}
\label{ghghgh} Let $F \in R[[q,...,q^{(r)}]]$ be a series such
that $G:=\Sigma_pF$ is $\d-$ symmetric. Then define the action of
the {\it Hecke operator} $T(p)_{\infty}$ on $F$ by
\begin{equation}
\label{suge} (T(p)_{\infty}F):= F(q^p,...,\d^r(q^p))+
G_{(p)}(0,...,0,q,...,0,...,0,q^{(r)}) \in R[[q,...,q^{(r)}]].
\end{equation}
Morally this should correspond to the Hecke action on
$\d-$expansions of weight (of degree) $0$.
\end{definition}

\begin{remark}
\label{rem1} If $G=\Sigma_pF$ is $\d-$symmetric then so is
$G^{\phi}=\Sigma_p(F^{\phi})$ and we have
\[(G^{\phi})_{(p)}=(G_{(p)})^{\phi}.\]
In particular $T(p)_{\infty}$ commutes with $\phi$ in the sense
that
\[T(p)_{\infty}(F^{\phi})=(T(p)_{\infty}F)^{\phi}\]
for any $F$ for which $\Sigma_pF$ is $\d-$symmetric.
\end{remark}

\begin{remark}
\label{rem2} If $F=\sum_{n \geq 1} c_nq^n \in R[[q]]$ then
$\Sigma_pF$ is $\d-$symmetric and
\begin{equation}
\label{baga} T(p)_{\infty}F=\sum_{n \geq 1} c_nq^{np}+p\sum_{n
\geq 1} c_{np}q^n.
\end{equation}
Indeed note that if
\begin{equation}
\label{oouf} q_1^n+...+q_p^n=P_n(S_1,...,S_p) \end{equation}
 with
$P_n$ a weighted homogeneous polynomial with $\bZ-$coefficients of
degree $n$ (with respect to the weights $1,2,...,p$) then
$P_n(0,...,0,q)$ is either $mq^{n/p}$ (with $m \in \bZ$) or $0$
according as $p$ divides $n$ or not. In case $n/p \in \bZ$,
specializing $q_i \mapsto \zeta_p^i$, $\zeta_p:=e^{2 \pi i/p}$, we
get $S_i \mapsto 0$ for $0 \leq i \leq p-1$ and $S_p \mapsto 1$ so
Equation \ref{oouf} yields $p=m$. Equation \ref{baga} follows.
Formula \ref{baga} is, in some sense, what one would  expect
 the action of $T(p)$ to yield
on series of order $0$ and weight $0$.
\end{remark}

\begin{remark}
\label{rem3} One can introduce a variant over $K$ of the above
definitions. A series
\[G \in K[[q_1,...,q_m,...,q_1^{(r)},...,q_m^{(r)}]]\]
is {\it $K-\d-$symmetric} if there exists a series
\[G_{(m)} \in K[[q_1,...,q_m,...,q_1^{(r)},...,q_m^{(r)}]]\]
such that Equation \ref{rupe} holds. (The series $G_{(m)}$ is,
again,  trivially seen to be unique; cf. \cite{dcc}.) For any
series $F \in K[[q,...,q^{(r)}]]$ one can define the series
$\Sigma_mF \in K[[q_1,...,q_m,...,q_1^{(r)},...,q_m^{(r)}]]$ as
before. If $F$ is such that $\Sigma_pF$ is $K-\d-$symmetric we
define $T(p)_{\infty}F \in K[[q,...,q^{(r)}]]$ by the formula
\ref{suge}. Then Remarks \ref{rem1} and \ref{rem2} hold verbatim
with $R$ replaced by $K$ and the words ``$\d-$symmetric" replaced
by ``$K-\d-$symmetric".
\end{remark}

\begin{definition}
A series $F \in R[[q,...,q^{(r)}]]$ is an {\it eigenvector} for
$T(p)_{\infty}$ with {\it eigenvalue} $\lambda \in R$ if
$\Sigma_pF$ is $\d-$symmetric and
\begin{equation}
\label{putza} T(p)_{\infty}F=\lambda \cdot F.
\end{equation}
\end{definition}

\begin{proposition}
\label{atpp}
 Let $f=\sum a_nq^n$ and $f^{\sharp}$ be as in
Theorems \ref{fur} and \ref{fur2} respectively. Then
$f^{\sharp}_{\infty,\rho}$ is an eigenvector of $T(p)_{\infty}$
with eigenvalue $a_p=a_p(f)$.
\end{proposition}

\begin{proof}
By Corollary \ref{ddus} $\Sigma_pf^{\sharp}_{\infty,\rho}$ is
$\d-$symmetric. Now we think of $f^{\sharp}_{\infty,\rho}$ as an
element of $K[[q,q',q'']]$ or $K[[q,q']]$ respectively and we
consider the extension of $T(p)_{\infty}$ ``over $K$" discussed in
Remark \ref{rem3}. Since $T(p)_{\infty}$ commutes with $\phi$ it
is enough to check that
\[\sum \frac{a_n}{n} q^n\]
is an eigenvector of $T(p)_{\infty}$ with eigenvalue $a_p$. This
can be checked directly as follows. First note that, by Equations
\ref{tzu} we have
\[pa_{n/p}+a_{np}=a_pa_n,\ \ \ n \geq 1.\]
Then, by Remark \ref{rem2}, we have
\[\begin{array}{rcl}
T(p)_{\infty} \left( \sum \frac{a_n}{n} q^n \right) & = &  \sum
\frac{a_n}{n} q^{np}+p \sum \frac{a_{np}}{np}q^n\\
\  & \  & \  \\
\  & = & \sum \frac{a_{n/p}}{n/p} q^n+\sum \frac{a_{np}}{n}q^n\\
\  & \  & \  \\
\ & = & \sum \frac{pa_{n/p}+a_{np}}{n}q^n\\
\  & \  & \  \\
\  & = &  a_p \sum\frac{a_n}{n}q^n.\end{array}\]
\end{proof}

One can develop a variant of $T(p)_{\infty}$ by allowing it to act
on certain series with denominators. We need a preparation. We
continue to denote by $S_1,...,S_p$  the fundamental symmetric
polynomials in $q_1,...,q_p$ and we let $s_1,...,s_p$ be
variables.

\begin{lemma}
\label{megga} Consider the $R-$algebras
\[
\begin{array}{rcl}
A & := &
R[[s_1,...,s_p]][s_p^{-1}]\h[s_1',...,s_p',...,s_1^{(r)},...,s_p^{(r)}]\h,\\
B & := &
R[[q_1,...,q_p]][q_1^{-1}...q_p^{-1}]\h[q_1',...,q_p',...,q_1^{(r)},...,q_p^{(r)}]\h.
\end{array}\] Then the natural algebra map \[A \ra B,\ \ \ s^{(i)}_j
\mapsto \d^iS_j\]
 is injective with torsion free cokernel.
\end{lemma}

\begin{proof}
Let $\sigma_1,...,\sigma_{p-1}$ be variables and let
$\Sigma_1,...,\Sigma_{p-1} \in R[S_1,...,S_p]$ be defined by
\[\Sigma_i:=q_1^i+...+q_p^i.\]
Note that
\[R[\Sigma_1,...,\Sigma_{p-1},S_p]=R[S_1,...,S_p].\]
So there is a natural isomorphism $C \simeq A$ where
\[C:=R[[\sigma_1,...,\sigma_{p-1},s_p]][s_p^{-1}]\h[\sigma_1',...,\sigma'_{p-1},
s_p',...,\sigma_1^{(r)},...,\sigma_{p-1}^{(r)},s_p^{(r)}]\h\] such
that the composition $C \ra A \ra B$ is given by
\[\sigma_j^{(i)} \mapsto \d^i \Sigma_j,\ \ s_p^{(i)} \mapsto \d^i
S_p.\] So it is enough to prove that $C \otimes k \ra B \otimes k$
is injective. We have
\[\begin{array}{rcl}
C \otimes k & = &
k[[\sigma_1,...,\sigma_{p-1},s_p]][s_p^{-1}][\sigma_1',...,\sigma'_{p-1},
s_p',...,\sigma_1^{(r)},...,\sigma_{p-1}^{(r)},s_p^{(r)}],\\
\  & \  & \  \\
B \otimes k & = & k[[q_1,...,q_p]][q_1^{-1}...q_p^{-1}]
[q_1',...,q_p',...,q_1^{(r)},...,q_p^{(r)}].\end{array}\] Now the
morphism
\[k[\sigma_1,...,\sigma_{p-1},s_p]=k[s_1,...,s_p] \ra
k[q_1,...,q_p]\] is finite and flat, and $(q_1,...,q_p)$ is the
unique maximal ideal lying over
\[(\sigma_1,...,\sigma_{p-1},s_p).\]
 Hence the ring homomorphism
\[k[[\sigma_1,...,\sigma_{p-1},s_p]] \ra
k[[q_1,...,q_p]]\] is faithfully flat, hence injective, so we have
an inclusion $L \subset M$ of their fraction fields. It is enough
to show that the map
\[L[\sigma_1',...,\sigma'_{p-1},
s_p',...,\sigma_1^{(r)},...,\sigma_{p-1}^{(r)},s_p^{(r)}] \ra
M[q_1',...,q'_p,...,q_1^{(r)},...,q_p^{(r)}]\] in injective. We
will show (and this will end our proof) that for each $i=0,...,r$
the images of
\begin{equation}
\label{ueito} \d^i \Sigma_1,...,\d^i \Sigma_{p-1},\d^i S_p \in
R[q_1,...,q_p,...,q_1^{(i)},...,q_p^{(i)}]
\end{equation}
in the ring
\[M[q_1',...,q'_p,...,q_1^{(i)},...,q_p^{(i)}]\]
are algebraically independent over
\[M[q_1',...,q'_p,...,q_1^{(i-1)},...,q_p^{(i-1)}].\]
Now one checks by induction on $i$ that for all $a=1,...,p-1$,
\[\d^i \Sigma_a=\sum_{j=1}^p (aq_j^{a(p-1)})^{p^i}
q_j^{(i)}+O(i-1)+pO(i),\] where
\[O(i) \in R[q_1,...,q_p,...,q_1^{(i)},...,q_p^{(i)}],\]
and $O(i-1)$ has the corresponding meaning. Similarly one has
\[\d^i S_p=\sum_{j=1}^p q_j^{(i)} (s_p/q_j)^{p^i}
+pO(i)+pO(i-1).\] So the images of the polynomials \ref{ueito} in
the ring \[k[q_1,...,q_p,...,q_1^{(i)},...,q_p^{(i)}]\] are
(non-homogeneous) linear polynomials in $q_1^{(i)},...,q_p^{(i)}$
with coefficients in
\[k[q_1,...,q_p,...,q_1^{(i-1)},...,q_p^{(i-1)}].\]
So we need to check that the matrix of the coefficients of
$q_1^{(i)},...,q_p^{(i)}$ in the reductions mod $p$
 of the
polynomials \ref{ueito} is non-singular. But this matrix is the
$p^i-$th power of the matrix
\[\left(
\begin{array}{ccccc}
1 & \cdot & \cdot & \cdot & 1\\
2q_1^p & \cdot & \cdot & \cdot & 2q_p^p\\
\cdot & \cdot & \cdot & \cdot & \cdot\\
\cdot & \cdot & \cdot & \cdot & \cdot\\
(p-1)q_1^{(p-2)p} & \cdot & \cdot & \cdot & (p-1)q_p^{(p-2)p}\\
s_p/q_1 & \cdot & \cdot & \cdot & s_p/q_p
\end{array}
\right)
\]
which is clearly non-singular.
\end{proof}

\begin{definition}
In the notations of Lemma \ref{megga}, an element  $G \in B$ will
be called {\it Laurent $\d-$symmetric} if it is the image of some
element $G_{(p)} \in A$ (which is then unique by  Lemma
\ref{megga}). For any $F \in R((q))\h[q',...,q^{(r)}]\h$ such that
\[\Sigma_pF:=\sum_{j=1}^pF(q_j,...,q_j^{(r)}) \in B\]
is Laurent $\d-$symmetric we may define
\[
 T(p)_{\infty,m}F:=
F(q^p,...,\d^r(q^p))+ p^mG_{(p)}(0,...,0,q,...,0,...,0,q^{(r)})
\in R[[q,...,q^{(r)}]].\] We write
$T(p)_{\infty}F=T(p)_{\infty,0}F$.
 A series $F \in
R((q))\h[q',...,q^{(r)}]\h$ is a {\it Laurent eigenvector} of
$T(p)_{\infty,m}$ with {\it eigenvalue} $\lambda \in R$ if
$\Sigma_pF$ is Laurent $\d-$symmetric and $T(p)_{\infty,m}F=
\lambda \cdot F$. The various values of $m$ morally correspond to
the Hecke action on $\d-$Fourier expansions of $\d-$modular forms
of various weights $w$ (with $deg(w)=-m$).
\end{definition}

\begin{remark}
If $F \in R[[q]][q',...,q^{(r)}]\h$ then one can consider the
following conditions:

a) $\Sigma_pF$ is $\d-$symmetric;

b) $\Sigma_pF$ is Laurent $\d-$symmetric.

A priori none of these conditions seems to imply the other. On the
other hand one can trivially see that $f^{\sharp}$ in Theorems
\ref{fur} and \ref{fur2} is not only $\d-$symmetric but also
Laurent $\d-$symmetric and a Laurent eigenvector for
$T(p)_{\infty}$ with eigenvalue $a_p(f)$.
\end{remark}

\begin{proposition}
For any $r \geq 1$ the $\d-$Fourier expansion
\[f^r_{\infty,\rho} \in R((q))\h[q',...,q^{(r)}]\h\]
of the $\d-$modular form
\[f^r=f^r_{crys} \in M^r(\Gamma_1(N),R,0)\]
is a Laurent eigenvector for $T(p)_{\infty,2}$ with eigenvalue
$p(p+1)$.
\end{proposition}

\begin{proof}
By Lemma \ref{olema} it is enough to show that $\Psi^{\phi^i}$ are
Laurent eigenvectors for $T(p)_{\infty,2}$ with eigenvalues
$p(p+1)$. Now $\Psi^{\phi^i}$ is Laurent $\d-$symmetric because
\[\begin{array}{rcl}
\sum_{j=1}^p \Psi^{\phi^i}(q_j,...,q_j^{(i+1)}) & = & \phi^i\left(
\sum_{j=1}^p \frac{1}{p} log\left(1+p \frac{q'_j}{q_j^p} \right)
\right)\\
\  & \  & \  \\
\ & = & \phi^i \left( \sum_{j=1}^p \frac{1}{p} log
\frac{\phi(q_j)}{q_j^p} \right)\\
\  & \  & \  \\
\  & = &  \phi^i \left( \frac{1}{p} log \prod_{j=1}^p
\frac{\phi(q_j)}{q_j^p} \right)\\
\  & \  & \  \\
\  & = & \phi^i \left( \frac{1}{p} log \frac{\phi(s_p)}{s_p^p}
\right)\\
\  & \  & \  \\
\  & = & \Psi^{\phi^i}(s_p).\end{array}\] Moreover, since
\[\Psi^{\phi^i}(q^p,...,\d^{i+1}(q^p))=
\phi^i\left( \frac{1}{p} log \frac{\phi(q^p)}{q^{p^2}} \right)= p
\cdot \Psi^{\phi^i},\] we have
\[T(p)_{\infty,2} \Psi^{\phi^i}=p \cdot \Psi^{\phi^i}+p^2 \cdot
\Psi^{\phi^i} =p(p+1)\Psi^{\phi^i}.\]
\end{proof}

We end by stating a result (to be proved in a subsequent paper)
showing that   there is an interesting relationship between
(Laurent) $\d-$symmetry and $\d-$ characters. This result is, in
some sense, a ``converse'' of the existence results for
$f^{\sharp}$ in the present paper. We assume, in what follows,
that $X_1(N)$ has genus at least $2$.

\begin{theorem}
\label{rem4} Fix an embedding $\rho:\bZ[1/N,\zeta_N] \ra R_p$ and
a $\d-$modular form $G \in M^r(\Gamma_1(N),R,0)$. Assume $G$
 is $\d-$holomorphic, $G$ vanishes at $\infty$,
and $\Sigma_pG_{\infty,\rho}$ is either $\d-$symmetric or Laurent
$\d-$symmetric. Then $G=\psi \circ \beta$ where $\beta:X_1(N)_R
\ra J_1(N)_R$ is the Abel-Jacobi map (corresponding to $\infty$)
and $\psi:J_1(N)_R \ra {\bf G}_{a,R}$ is a $\d-$character. (In
particular $G$ is automatically a $\d-$cusp form.)
\end{theorem}

The case when $\Sigma_pG_{\infty,\rho}$ is $\d-$symmetric follows
directly from the main Theorem of \cite{dcc}. The case when
$\Sigma_pG_{\infty,\rho}$ is Laurent $\d-$symmetric can be proved
in an entirely similar way.

\begin{remark}
Given a classical newform $f=\sum a_n q^n \in
S_2(\Gamma_0(N),\bC)$, an embedding $\rho:\cO_F[1/MN,\zeta_N] \ra
R_p$ (where the cusps are defined over $\cO_F[1/M]$), and an
embedding $\e:\cO_f^{(Np)} \ra R_p$ one is naturally lead to try
to compute the $R-$module ${\mathcal M}={\mathcal M}(f,\rho,\e,
r,p)$ of all $\d-$ modular forms $G \in M^r(\Gamma_1(N),R_p,0)$
satisfying the following properties:

1) $G$ is a $\d-$cusp form (with respect to $\rho$),

2) $G$ belongs to $f$ (outside $Np$) with character $\e$ and
exponent $0$,

3) $G_{\infty,\rho}$ is an eigenvector (respectively a Laurent
eigenvector) of $T(p)_{\infty}$ with eigenvalue $\e(a_p)$.

By Proposition \ref{rem4} and Theorem \ref{mythm}, $rank\
{\mathcal M} \leq rg_1(N)$ where $g_1(N)$ is the genus of
$X_1(N)_{\bf C}$. One should expect a much better bound for $rank\
B$. Of course, by Theorems \ref{fur} and \ref{fur2} plus
Proposition \ref{atpp}, if $g=[K_f:\bQ]=1$ then
$f^{\sharp}_{\infty,\rho} \in {\mathcal M}$.
\end{remark}

\bibliographystyle{amsplain}


\end{document}